\documentclass[aos,MSNbibl,seceqn,dvips]{arximspdf}

\doi{10.1214/12-AOS1073} 
\volume{41}
\issue{1}
\pubyear{2013}
\firstpage{269}
\lastpage{295}

\makeatletter

\newcommand{\cal}{\mathcal}

\newtheorem{theorem}{Theorem}[section]
\newtheorem{lemma}{Lemma}[section]
\newtheorem{cor}{Corollary}[section]

\newproclaim{remark}{Remark}[section]
\newproclaim{cond}{Condition}[section]

\makeatother

\begin{document}
\begin{frontmatter}

\title{Weighted likelihood estimation under two-phase~sampling}
\runtitle{Weighted likelihood: Two-phase sampling}

\begin{aug}
\author[A]{\fnms{Takumi} \snm{Saegusa}\corref{}\thanksref{t1}\ead[label=e1]{tsaegusa@uw.edu}}
\and
\author[B]{\fnms{Jon A.} \snm{Wellner}\thanksref{t2}\ead[label=e2]{jaw@stat.washington.edu}}
\runauthor{T. Saegusa and J. A. Wellner}
\affiliation{University of Washington}
\address[A]{Department of Biostatistics\\
University of Washington\\
Seattle, Washington 98195-7232\\
USA\\
\printead{e1}}
\address[B]{Department of Statistics\\
University of Washington\\
Seattle, Washington 98195-4322\\
USA\\
\printead{e2}} 
\end{aug}

\thankstext{t1}{Supported by NIH/NIAID R01 AI089341.}
\thankstext{t2}{Supported in part by NSF Grant DMS-11-04832,
NI-AID Grant 2R01 AI291968-04
and the Alexander von Humboldt Foundation.}

\received{\smonth{1} \syear{2012}}
\revised{\smonth{8} \syear{2012}}

%
\begin{abstract}
We develop asymptotic theory for weighted likelihood estimators (WLE)
under two-phase stratified sampling without replacement. We also
consider several variants of WLEs involving estimated weights and
calibration. A set of empirical process tools are developed including a
Glivenko--Cantelli theorem, a theorem for rates of convergence of
$M$-estimators, and a Donsker theorem for the inverse probability
weighted empirical processes under two-phase sampling and sampling
without replacement at the second phase. Using these general results,
we derive asymptotic distributions of the WLE of a finite-dimensional
parameter in a general semiparametric model where an estimator of a
nuisance parameter is estimable either at regular or nonregular rates.
We illustrate these results and methods in the Cox model with right
censoring and interval censoring. We compare the methods via their
asymptotic variances under both sampling without replacement and the
more usual (and easier to analyze) assumption of Bernoulli sampling at
the second phase.
\end{abstract}

%

\end{frontmatter}

\section{Introduction}\label{sec1}
Two-phase sampling is a sampling technique that aims at cost reduction
and improved efficiency of estimation.
At phase I, a large sample is drawn from a population, and information on
variables that are easier to measure is collected.
These phase I variables may be important variables such as exposure in
a regression model,
or simply may be
auxiliary variables that are correlated with unavailable variables at
phase I.
The sample space is then stratified based on these phase I variables.
At phase II, a subsample is drawn without replacement from each stratum
to obtain phase II
variables that are costly or difficult to measure.
Strata formation seeks either to oversample subjects with important
phase I variables,
or to effectively sample subjects with targeted phase II variables, or both.
This way, two-phase sampling achieves effective access to important variables
with less cost.

While two-phase sampling was originally introduced
in survey sampling
by Neyman
\cite{Neyman38} for estimation of the ``finite population mean'' of
some variable,
it has become increasingly important in a variety of areas of
statistics, biostatistics and epidemiology,
especially since~\cite{White1982,Prentice86} and~\cite{MR924857}.

The setting treated here is as follows:

$\bullet$
We begin with a semiparametric model ${\cal P}$ for a vector of
variables $X$ with values in ${\cal X}$.
[The prime examples which we treat in detail in Section~\ref{sec4} are
the Cox
proportional hazards
regression model with (a) right censoring, and (b) interval
censoring.]

$\bullet$ 
Let $W = (X,U) \in{\cal X} \times{\cal U} \equiv{\cal W}$ where $U$
is a vector of
``auxiliary variables,'' not involved in the model ${\cal P}$.
Suppose that $W \sim\tilde{P}_0$ and $X \sim P_0$.
Now suppose that $V \equiv(\tilde{X}, U) \in{\cal V}$ where $\tilde
{X} \equiv\tilde{X} (X)$ is a
coarsening of $X$.

$\bullet$ 
At phase I we observe $V_1,\ldots, V_N $ i.i.d. as $V$, and then use
the phase I
data to form strata, that is, disjoint subsets ${\cal V}_1,\ldots,
{\cal V}_J$ of ${\cal V}$ with
$\sum_{j=1}^J {\cal V}_j = {\cal V}$.
We let $N_j = \# \{ i \le N\dvtx V_i \in{\cal V}_j \}$.

$\bullet$ 
Next, a phase II sample is drawn by sampling without replacement $n_j
\le N_j$ items from stratum $j$.
For the items selected we observe $X_i$.
Thus for the selection indicators $\xi_i$ we have
$\tilde{P}_0 ( \xi_i = 1 | V_i ) = (n_j/N_j) 1_{{\cal V}_j} (V_i )
\equiv\pi_0 (V_i)$.


$\bullet$ 
Finally weighted likelihood (or inverse probability weighted)
estimation methods based on all the observed data are used to estimate
the parameters of the model ${\cal P}$ and to make further inferences
about the model.

It is now well known that the classical Horvitz--Thompson estimators~\cite{MR0053460}
use only the phase II data and are inefficient,
sometimes quite severely so; see, for example,
\cite{MR1294730,Lumley2011,Breslow2009a,Breslow2009b} and
\cite{Zheng-Little04}. Improvements in efficiency of estimation can be
achieved by ``adjusting'' the weights by use of the phase I data (even
though the sampling probabilities are known). Two basic methods of
adjustment are:

(1) Estimated weights, a method originating in the missing data literature
\cite{MR1294730}, and with significant further developments since in connection
with many models in which the missing-ness mechanism is not known, in
contrast to our
current setting in which the missing-ness is by design.

(2) Calibration, a method originating in the sample survey literature
\cite{MR1173804}; see also~\cite{LumleyRbook,Lumley2011}.

One of our goals here is to study existing methods for adjustment of
the weights of weighted likelihood
methods and to introduce several new methods: modified calibration as
suggested by Chan
\cite{chan12} 
and centered calibration as proposed here in Section~\ref{sec2}.

A second goal is to give a systematic treatment of estimators based on
sampling without replacement
at phase II in the setting of general semiparametric models
and to make comparisons with the behavior of estimators based on
Bernoulli (or independent)
sampling at phase II, thus continuing and strengthening the comparisons
made in
\cite{MR2325244,MR2391566},
and~\cite{Breslow2009a,Breslow2009b} for a particular sub-class of
semiparametric models
and adjustments via estimated weights and ordinary calibration.
Many studies of the theoretical properties of
procedures based on two-phase design data have been made for the case
of Bernoulli sampling;
see, for example,~\cite{MR2860827} and the review of case-cohort
sampling given there.
On the other hand, while statistical practice continues to involve
phase II data sampled without replacement,
most available theory in this case (other than \cite
{MR2325244,MR2391566}) has
developed on a model-by-model basis.
As has become clear from~\cite{MR2325244,MR2391566}, sampling without
replacement
results in smaller asymptotic variances, 
and hence inference based on asymptotic variances derived from
Bernoulli sampling will often be conservative.
Our treatment here provides theory and tools for dealing directly with
the sampling without replacement
design.
We do this by providing the relevant
theory both for semiparametric models in which the infinite-dimensional
nuisance parameters
can be estimated at a regular rate $(\sqrt{n})$ with complete data, and
semiparametric models
in which the infinite-dimensional nuisance parameters can only be
estimated at slower (nonregular) rates.

The main contributions of our paper are three-fold:
First, we establish two $Z$-theorems giving weak sufficient conditions
for asymptotic distributions of the WLEs in general semiparametric models.
The first theorem covers the case where the nuisance parameter is estimable
at a regular rate; this yields rigorous justification of
\cite{Breslow2009a,Breslow2009b} under weaker conditions.
The second theorem covers the case of general semiparametric models
with nonregular rates for estimators of the nuisance parameters.
The conditions of our theorems, formulated in terms of complete data,
are almost identical to those for the MLE with complete data.
This formulation allows us to use tools from empirical process theory together
with the new tools developed here in a straightforward way.
Second, we propose centered calibration, a new calibration method.
This new calibration method is the only one guaranteed
to yield improved efficiency over the plain WLE under both
Bernoulli sampling and sampling without replacement, while other
methods are
warranted only for Bernoulli sampling.
Third, we establish general results for the inverse probability weighted
(IPW) empirical process.
Some results such as a Glivenko--Cantelli theorem (Theorem~\ref{thmfpsPiGC})
and a Donsker theorem (Theorem~\ref{thmfpsD}) are of interest in their
own right.
These results accounting for dependence due to the sampling design are
useful in verifying the conditions of $Z$-theorems in applications.
For instance, Theorems~\ref{thmfpsPiGC} and~\ref{thmrate} easily
establish consistency and rates of convergence under our ``without
replacement'' sampling scheme.
We illustrate application of the general results with examples in
Section~\ref{sec4}.\vadjust{\goodbreak}

The rest of the paper is organized as follows.
In Section~\ref{sec2}, we introduce our 
estimation procedures
in the context of a general semiparametric model.
The WLE and methods involving adjusted weights are discussed.
Two $Z$-theorems are presented in Section~\ref{sec3}; these yield
asymptotic distributions
of the WLEs of finite-dimensional parameters of the model.
All estimators are compared under Bernoulli sampling and sampling
without replacement with different methods of adjusting weights.
In Section~\ref{sec4} we apply our $Z$-theorems to the Cox model with
both right
censoring and interval censoring.
Section~\ref{sec5} consists of general results for IPW empirical processes.
Several open problems are briefly discussed in Section~\ref{sec6}.
All proofs, except those in Section~\ref{sec4} and auxiliary results, are
collected in
\cite{Saegusa-Wellner2012supp}.

\section{Sampling, models and estimators}\label{sec2}
We use the basic notation introduced in the previous section.
After stratified sampling, $X$ is fully observed for $n_j$ subjects in
the $j$th stratum at phase II.
The observed data is $(V,X\xi,\xi)$ where $\xi$ is the indicator of
being sampled at phase II.
We use a doubly subscripted notation: for example, 
$V_{j,i}$ 
denotes $V$ for the $i$th subject in stratum~$j$.
We denote the stratum probability for the $j$th stratum by
$\nu_j\equiv\tilde{P}_0(V\in\mathcal{V}_j)$,
and the conditional expectation given membership in the $j$th stratum by
$P_{0|j}(\cdot) \equiv\tilde{P}_0(\cdot|V\in\mathcal{V}_j)$.

The sampling probability is $P(\xi=1|V_i)=\pi_0(V_i)=n_j/N_j$ for
$V_i\in\mathcal{V}_j$.
These sampling probabilities are assumed to be strictly positive;
that is, there is a
constant $\sigma>0$ such that
$0<\sigma\leq\pi_0(v)\leq1$ for $v\in\mathcal{V}$.
We assume that $n_j/N_j\rightarrow p_j>0$ for $j=1,\ldots, J$ as
$N\rightarrow\infty$.
Although dependence is induced among the observations $(V_i,\xi
_iX_i,\xi
_i)$ by the sampling
indicators, the vector of sampling indicators $(\xi_{j1},\ldots,\xi
_{jN_j})$ within strata,
are exchangeable for each $j=1,\ldots, J$, and the $J$ random vectors
$(\xi_{j1},\ldots,\xi_{jN_j})$ are independent.

The empirical measure is one of the most useful tools in empirical
process theory.
Because the $X_i$'s are observed only for a sub-sample at phase II,
we define, instead,
the IPW
empirical measure $\mathbb{P}_N^{\pi}$ by
\[
\mathbb{P}_N^\pi=\frac{1}{N}\sum
_{i=1}^N\frac{\xi_i}{\pi_0(V_i)}\delta_{X_i} =
\frac{1}{N}\sum_{j=1}^J\sum
_{i=1}^{N_j}\frac{\xi_{j,i}}{n_j/N_j}\delta_{X_{j,i}},
\]
where $\delta_{X_i}$ denotes a Dirac measure placing unit mass on $X_i$.
The identity in the last display is justified by the arguments in
Appendix A of~\cite{MR2325244}.
We also define the IPW empirical process by $\mathbb{G}_N^\pi=\sqrt
{N}(\mathbb{P}_{N}^{\pi}-P_0)$
and the phase II empirical process for the $j$th stratum by
$\mathbb{G}_{j,N_j}^\xi
\equiv\sqrt{N_j}(\mathbb{P}_{j,N_j}^\xi-(n_j/N_j)\mathbb{P}_{j,N_j})$,
$j=1,\ldots,J$,
$\mathbb{P}_{j,N_j}^\xi\equiv N_j^{-1}\sum_{i=1}^{N_j}\xi
_{j,i}\delta_{X_{j,i}}$
is the phase II empirical measure for the $j$th stratum, and
$\mathbb{P}_{j,N_j}\equiv N_j^{-1}\sum_{i=1}^{N_j}\delta_{X_{j,i}}$
is the empirical measure for all the data in the $j$th stratum; note
that the latter empirical measure
is not observed.
Then, following~\cite{MR2325244},
we decompose $\mathbb{G}_N^\pi$ as follows:
%
%
\begin{equation}\label{BasicDecompositionBreslowWellner}
\mathbb{G}_N^\pi= \mathbb{G}_N +\sum
_{j=1}^J \sqrt{\frac{N_j}{N}} \biggl(
\frac{N_j}{n_j} \biggr)\mathbb{G}_{j,N_j}^\xi,
\end{equation}
where $\mathbb{P}_N = N^{-1} \sum_{j=1}^J N_j \mathbb{P}_{j, N_j}$
and $\mathbb{G}_N = \sqrt{N} (\mathbb{P}_N - P_0)$.
Notice that
$\mathbb{G}_{j,N_j}^{\xi}$ correspond to
``exchangeably weighted bootstrap'' versions of the stratum-wise
complete data
empirical processes $\mathbb{G}_{j, N_j} \equiv\sqrt{N_j} ( \mathbb
{P}_{j, N_j} - P_{0|j})$.
This observation
allows application of the ``exchangeably weighted bootstrap'' theory of
\cite{MR1245301} and~\cite{MR1385671}, Section 3.6.

\subsection{Improving efficiency by adjusting weights}\label{sec2.1}
Efficiency of estimators based on IPW empirical processes can be
improved by adjusting weights,
either by estimated weights
\cite{MR1294730} or by calibration~\cite{MR1173804}
via use of the phase I data; see also~\cite{Lumley2011}.
Besides these, we discuss two variants of calibration,
modified calibration~\cite{chan12}, and our proposed new method,
centered calibration.

Let $Z_i\equiv g(V_i)$ be the auxiliary variables for the $i$th subject
for a known transformation $g$.
For estimated weights with binary regression,
$Z_i$ contains the membership indicators for the strata
$I_{\mathcal{V}_j}(V_i), j=1,\ldots,J$.
Observations with $\pi_0(V)=1$ are dropped from binary regression,
and the original weight $1$ is used.
For notational simplicity, we write $Z_i$ for either method, and assume
that sampling
probabilities are strictly less than 1 for all strata.

\subsubsection{Estimated weights}\label{sec2.1.1}
The method of estimated weights adjusts weights
through binary regression on the phase I variables.
The sampling\vspace*{1pt} probability for the $i$th subject is modeled by
$p_\alpha(\xi_i|Z_i)=G_e(Z_i^T\alpha)^{\xi_i}(1-G_e(Z_i^T\alpha
))^{1-\xi_i}
\equiv\pi_{\alpha}(V_i)^{\xi_i} \{1-\pi_\alpha(V_i) \}^{1-\xi_i}$,
where $\alpha\in\mathcal{A}_e\subset\mathbb{R}^{J+k}$ is a regression
parameter and
$G_e\dvtx\mathbb{R}\mapsto[0,1]$ is a known function.
If $G_e(x)=e^x/(1+e^x)$, for instance, then the adjustment simply
involves logistic regression.
Let $\hat{\alpha}_N$ be the estimator of $\alpha$ that maximizes the
pseudo- (or composite) likelihood
%
%
\begin{equation}
\label{eqncomplik} \prod_{i=1}^Np_\alpha(
\xi_i|Z_i) =\prod_{i=1}^NG_e
\bigl(Z_i^T\alpha\bigr)^{\xi_i}
\bigl(1-G_e\bigl(Z_i^T\alpha\bigr)
\bigr)^{1-\xi_i}.
\end{equation}
We define the IPW empirical measure with estimated weights by
\[
\mathbb{P}_{N}^{\pi,e} =\frac{1}{N}\sum
_{i=1}^N\frac{\xi_i}{\pi_{\hat{\alpha
}_N}(V_i)}\delta_{X_i} =
\frac{1}{N}\sum_{i=1}^N
\frac{\xi_i}{\pi_0(V_i)}\frac{\pi_0(V_i)}{G_e(Z^T_i\hat{\alpha
}_N)}\delta_{X_i},
\]
and the IPW empirical process with estimated weights by
$\mathbb{G}_{N}^{\pi,e}=\break\sqrt{N}(\mathbb{P}_{N}^{\pi,e}-P_0)$.

\subsubsection{Calibration}\label{sec2.1.2}
Calibration adjusts weights so that the inverse probability weighted average
from the phase\vadjust{\goodbreak} II sample is equated to the phase I average,
whereby the phase I information is taken into account for estimation.
Specifically, we find an estimator $\hat{\alpha}_N$ that is the
solution for
$\alpha\in\mathcal{A}_{c}\subset\mathbb{R}^k$ of the following
calibration equation:
%
%
\begin{equation}
\label{eqncaleqn2} \frac{1}{N}\sum_{i=1}^N
\frac{\xi_iG_{c}(V_i;\alpha)}{\pi_0(V_i)}Z_i =\frac{1}{N}\sum
_{i=1}^NZ_i,
\end{equation}
where $G_{c}(V;\alpha)\equiv G(g(V)^T\alpha)=G(Z^T\alpha)$ for known
$G$ with $G(0)=1$ and $\dot{G}(0)>0$.
We call $\pi_\alpha(V)\equiv\pi_0(V)/G_{c}(V;\alpha)$
the calibrated sampling probability.
We define the \textit{calibrated IPW empirical measure} by
\[
\mathbb{P}_{N}^{\pi,c} = \frac{1}{N}\sum
_{i=1}^N\frac{\xi_i}{\pi_{\hat{\alpha
}_N}(V_i)}\delta_{X_i} =
\frac{1}{N}\sum_{i=1}^N
\frac{\xi_i}{\pi_0(V_i)}G \bigl(Z^T_i\hat{\alpha
}_N \bigr)\delta_{X_i}
\]
and the \textit{calibrated IPW empirical process} by
$\mathbb{G}_{N}^{\pi,c}=\sqrt{N}(\mathbb{P}_{N}^{\pi,c}-P_0)$.

Examples 
for $G$ in the definition of $G_c$ are listed in~\cite{MR1173804} ($F$
in their notation).
For $G(x)=1+x$, $\mathbb{P}_N^{\pi,c} X$ is a well-known regression
estimator of 
the mean of
$X$.
Since we assume boundedness of $G$ later, we may want to
consider truncated versions of these examples instead.
Note that choice of $G$ in (variants of) calibration does not affect
asymptotic results on WLEs.

As noted in
\cite{LumleyRbook}, there are several different approaches to calibration.
Here, and in introducing variants of calibration below,
we adopt the view that
calibration proceeds by making 
the smallest
possible change in weights in order to match 
an estimated phase II
average with the corresponding
phase I average.
Another approach proceeds via
regression modeling of the variable $X$ of interest
and the auxiliary variables $V$, leading to a robustness discussion on
effects of the validity of the model on estimation for $X$.
We prefer the former view because we do not assume a model for
$X$ and $V$ throughout this paper.
In fact, our results are independent of such a modeling assumption.

\subsubsection{Modified calibration}\label{sec2.1.3}
\label{subsubsecimprovedcalibration}
Modifying the function $G_c$ in calibration so that individuals
with higher sampling probabilities $\pi(V_i)$ receive less weight
was proposed by
\cite{chan12} in a missing response problem where observations are i.i.d.
(see, e.g.,~\cite{MR2836413} for recent developments in this area and
\cite{Lumley2011} for their connections with calibration methods).
An interpretation of this method within the framework of~\cite{MR1173804}
is discussed in~\cite{Saegusa-Wellner2012TechRep}.
In modified calibration, we find the estimator $\hat{\alpha}_N$ that is
the solution for
$\alpha\in\mathcal{A}_{\mathrm{mc}}\subset\mathbb{R}^k$ of the following
calibration equation:
%
%
\begin{equation}
\label{eqnicaleqn} \frac{1}{N}\sum_{i=1}^N
\frac{\xi_iG_{\mathrm{mc}}(V_i;\alpha)}{\pi_0(V_i)}Z_i =\frac{1}{N}\sum
_{i=1}^NZ_i,
\end{equation}
where
$G_{\mathrm{mc}}(V;\alpha)\equiv
G((\pi_0(V)^{-1}-1)Z^T\alpha)$ for known $G$ with $G(0)=1$ and $\dot
{G}(0)>0$.
We call $\pi_\alpha(V)\equiv\pi_0(V)/G_{\mathrm{mc}}(V;\alpha)$ the \textit
{calibrated\vadjust{\goodbreak}
sampling probability with modified calibration}.
We define the \textit{IPW empirical measure with modified calibration} by
\[
\mathbb{P}_{N}^{\pi,\mathrm{mc}} = \frac{1}{N}\sum
_{i=1}^N\frac{\xi_i}{\pi_{\hat{\alpha
}_N}(V_i)}\delta_{X_i} =
\frac{1}{N}\sum_{i=1}^N
\frac{\xi_i}{\pi_0(V_i)}G \biggl(\frac
{1-\pi_0(V_i)}{\pi_0(V_i)}Z^T_i\hat{
\alpha}_N \biggr)\delta_{X_i}
\]
and the corresponding IPW empirical process by
$\mathbb{G}_{N}^{\pi,\mathrm{mc}}=\sqrt{N}(\mathbb{P}_{N}^{\pi,\mathrm{mc}}-P_0)$.

\subsubsection{Centered calibration}\label{sec2.1.4}
We propose a new method, centered calibration, that calibrates
on centered auxiliary variables with modified calibration.
This method improves the plain WLE under our sampling scheme,
while retaining the good properties of modified calibration.
See Section~\ref{sec3.4} for a discussion of its advantage and
connections to
other methods.

In centered calibration, we find the estimator $\hat{\alpha}_N$ that is
the solution for
$\alpha\in\mathcal{A}_{\mathrm{cc}}\subset\mathbb{R}^k$ of the following
calibration equation:
%
%
\begin{equation}
\label{eqnccaleqn} \frac{1}{N}\sum_{i=1}^N
\frac{\xi_iG_{\mathrm{cc}}(V_i;\alpha)}{\pi_0(V_i)}(Z_i-\overline{Z}_N) =0,
\end{equation}
where
$G_{\mathrm{cc}}(V;\alpha)\equiv G((\pi_0(V)^{-1}-1)(Z-\overline
{Z}_N)^T\alpha)$
for known $G$ with $G(0)=1$ and $\dot{G}(0)>0$ and $\overline
{Z}_N=N^{-1}\sum_{i=1}^NZ_i$.
We call $\pi_\alpha(V)\equiv\pi_0(V)/G_{\mathrm{cc}}(V;\alpha)$ the \textit
{calibrated
sampling probability with centered calibration}.
We define the \textit{IPW empirical measure with centered calibration} by
\[
\mathbb{P}_{N}^{\pi,\mathrm{cc}} = \frac{1}{N}\sum
_{i=1}^N\frac{\xi_i}{\pi_{\hat{\alpha
}_N}(V_i)}\delta_{X_i} =
\frac{1}{N}\sum_{i=1}^N
\frac{\xi_i}{\pi_0(V_i)}G_{\mathrm{cc}}(V_i;\hat{\alpha}_N)
\delta_{X_i}
\]
and the corresponding IPW empirical process by
$\mathbb{G}_{N}^{\pi,\mathrm{cc}}=\sqrt{N}(\mathbb{P}_{N}^{\pi,\mathrm{cc}}-P_0)$.

\subsection{Estimators for a semiparametric model ${\cal P}$}\label{sec2.2}

We study the asymptotic distribution of the weighted likelihood
estimator of a
finite-dimensional parameter $\theta$ in a general semiparametric model
$\mathcal{P}=\{P_{\theta,\eta}\dvtx\theta\in\Theta,\eta\in H\}$ where
$\Theta\subset\mathbb{R}^p$
and the nuisance parameter space $H$ is a subset of some Banach space
$\mathcal{B}$.
Let $P_0=P_{\theta_0,\eta_0}$ denote the true distribution.

The MLE for complete data is often obtained
as a solution to the infinite-dimensional likelihood equations.
In such models, the WLE under two-phase sampling
is obtained by solving the corresponding infinite-dimensional inverse
probability weighted likelihood equations.
Specifically, the WLE
$(\hat{\theta}_N,\hat{\eta}_N)$ is a solution to the following weighted
likelihood equations:
%
%
\begin{eqnarray}
\label{eqnwlik}
\Psi_{N,1}^\pi(\theta,\eta)&=&
\mathbb{P}_N^\pi\dot{\ell}_{\theta,\eta
}=o_{P^*}
\bigl(N^{-1/2}\bigr),
\nonumber\\[-8pt]\\[-8pt]
\bigl\|\Psi_{N,2}^\pi(\theta,\eta)h\bigr\|_{\mathcal{H}}
&=&\bigl\|\mathbb{P}_N^\pi(B_{\theta,\eta}h-P_{\theta,\eta
}B_{\theta,\eta
}h)
\bigr\|_{\mathcal{H}} =o_{P^*}\bigl(N^{-1/2}\bigr),\nonumber
\end{eqnarray}
where $\dot{\ell}_{\theta,\eta}\in\mathcal{L}_2^0(P_{\theta,\eta})^p$
is the score function for $\theta$, and the score operator
$B_{\theta,\eta}\dvtx\mathcal{H}\mapsto\mathcal{L}_2^0(P_{\theta,\eta})$
is the bounded linear\vadjust{\goodbreak} operator mapping a direction $h$ in some Hilbert
space $\mathcal{H}$ of one-dimensional submodels for $\eta$ along which
$\eta\rightarrow\eta_0$. The WLE with estimated weights
$(\hat{\theta}_{N,e},\hat{\eta}_{N,e})$, the calibrated WLE
$(\hat{\theta}_{N,c},\hat{\eta}_{N,c})$, the WLE with modified
calibration $(\hat{\theta}_{N,\mathrm{mc}},\hat{\eta}_{N,\mathrm{mc}})$
and the WLE with centered calibration
$(\hat{\theta}_{N,\mathrm{cc}},\hat{\eta}_{N,\mathrm{cc}})$ are
obtained by replacing $\mathbb{P}_N^\pi$ by $\mathbb{P}_N^{\pi,\#}$
with $\#\in\{ e,c,\mathrm{mc},\mathrm{cc}\}$ in (\ref{eqnwlik}),
respectively. Let $\dot{\ell}_0=\dot{\ell}_{\theta_0,\eta_0}$ and
$B_0=B_{\theta _0,\eta_0}$.

\section{Asymptotics for the WLE in general semiparametric
models}\label{sec3}
We consider two cases:
in the first case the nuisance parameter $\eta$ is estimable
at a regular (i.e., $\sqrt{n}$) rate, and
for ease of exposition, $\eta$ is assumed to be a measure.
In the second case $\eta$ is only estimable at a nonregular (slower
than $\sqrt{n}$) rate.
Our theorem (Theorem~\ref{thmzthm2}) concerning the second case nearly
covers the
former case, but requires slightly more smoothness and a separate proof of
the rate of convergence for an estimator of $\eta$. 
On the other hand, our theorem (Theorem~\ref{thmzthm1}) concerning the
former case includes a proof of the (regular) $(\sqrt{n})$ rate of convergence,
and hence is of interest by itself.

\subsection{Conditions for adjusting weights}\label{sec3.1}
We assume the following conditions for estimators of $\alpha$
for 
adjusted weights.
Throughout this paper, we may assume both Conditions~\ref{condest}
and~\ref{condcal} at the same time, but it should be understood that the
former condition is used exclusively for the estimators regarding estimated
weights and the latter condition is imposed only for estimators
regarding (variants of) calibration.
Also, it should be understood that Conditions~\ref{condcal}(a)(i) and (d)(i),
Conditions~\ref{condcal}(a)(ii)
and (d)(ii) and Conditions~\ref{condcal}(a)(iii) and (d)(iii) are assumed
for estimators defined via 
calibration, modified calibration and centered calibration, respectively.
%
%
\begin{cond}[(Estimated weights)]
\label{condest}
(a) The estimator $\hat{\alpha}_N$ is a maximizer
of the pseudo-likelihood (\ref{eqncomplik}).

(b) $Z\in\mathbb{R}^{J+k}$ is not concentrated on a
$(J+k)$-dimensional
affine space of $\mathbb{R}^{J+k}$ and has bounded support.

(c) $G_e\dvtx\mathbb{R}\mapsto[0,1]$ is a twice continuously
differentiable, monotone function.

(d)
$S_0\equiv P_0[\{\dot{G}_e(Z^T\alpha_0)\}^2\{\pi_0(V)(1-\pi_0(V))\}
^{-1}Z^{\otimes2}]$\vspace*{1pt}
is finite and nonsingular, where $\dot{G}_e$ is a derivative of $G_e$.

(e) The ``true'' parameter
$\alpha_0=(\alpha_{0,1},\ldots,\alpha_{0,J+k})$ is given by
$\alpha_{0,j}=G_e^{-1}(p_j)$, for $j=1,\ldots,J$ and $\alpha_{0,j}=0$,
for $j=J+1,\ldots,J+k$. The parameter $\alpha$ is identifiable, that
is, $p_\alpha= p_{\alpha_0}$ almost surely implies $\alpha= \alpha_0$.

(f) For a fixed $p_j \in(0,1)$, $n_j$ satisfies
$n_j=[N_jp_j]$ for $j=1,\ldots,J$.
\end{cond}

%
\begin{cond}[(Calibrations)]
\label{condcal}
(a) (i) The estimator $\hat{\alpha}_N = \hat{\alpha}_N^c$ is a
solution of calibration equation (\ref{eqncaleqn2}).
(ii) The estimator $\hat{\alpha}_N = \hat{\alpha}_N^{\mathrm{mc}}$ is a
solution of calibration equation (\ref{eqnicaleqn}).
(iii) The estimator $\hat{\alpha}_N = \hat{\alpha}_N^{\mathrm{cc}}$ is a
solution of calibration equation~(\ref{eqnccaleqn}).

(b) $Z\in\mathbb{R}^k$ is not concentrated at $0$ and has
bounded support.\vadjust{\goodbreak}

(c) $G$ is a strictly increasing continuously differentiable
function on
$\mathbb{R}$ such that $G(0)=1$ and for all $x$, $-\infty<m_1\leq G(x)
\leq M_1<\infty$
and $0< \dot{G}(x) \leq M_2<\infty$, where $\dot{G}$ is the derivative
of $G$.

(d) (i) $P_0Z^{\otimes2}$ is finite and positive definite.
(ii) $P_0[\pi_0(V)^{-1} (1-\pi_0(V))Z^{\otimes2}]$ is finite and
positive definite.
(iii) $P_0[\pi_0(V)^{-1} (1-\pi_0(V))(Z-\mu_Z)^{\otimes2}]$ is finite
and positive definite where $\mu_Z=PZ$.

(e) The ``true'' parameter $\alpha_0=0$.
\end{cond}

Condition~\ref{condest}(f)
may seem unnatural at first, but in practice the
phase II sample size $n_j$ can be chosen by the investigator
so that the sampling probability $p_j$ can be understood to be
automatically chosen to satisfy $n_j=[N_jp_j]$. The other parts of
Condition~\ref{condest}
are standard in binary regression, and Condition~\ref{condcal}
is similar to Condition~\ref{condest}.

Asymptotic properties of $\hat{\alpha}_N$ for all methods are
proved in~\cite{Saegusa-Wellner2012supp}.

\subsection{Regular rate for a nuisance parameter}\label{sec3.2}
We assume the following conditions.
%
%
\begin{cond}[(Consistency)]
\label{condwlereg1}
The estimator $(\hat{\theta}_N,\hat{\eta}_N)$ is consistent for
$(\theta_0,\eta_0)$ and
solves the weighted likelihood equations (\ref{eqnwlik}), where
$\mathbb{P}_N^\pi$
may be replaced by $\mathbb{P}_N^{\pi,\#}$ with $\#\in\{e,c,\mathrm{mc},\mathrm{cc}\}$
for the estimators with adjusted weights.
\end{cond}

%
\begin{cond}[(Asymptotic equicontinuity)]
\label{condwlereg2}
Let $\mathcal{F}_1(\delta)=\{\dot{\ell}_{\theta,\eta}\dvtx
|\theta
-\theta_0|
+\|\eta-\eta_0\|<\delta\}$
and
\mbox{$\mathcal{F}_2(\delta)=\{B_{\theta,\eta}h-P_{\theta,\eta
}B_{\theta,\eta
}h\dvtx h\in\mathcal{H},|\theta-\theta_0|
+\|\eta-\eta_0\|<\delta\}$}.
There exists a $\delta_0>0$ such that
(1) $\mathcal{F}_k(\delta_0),k=1,2$, are $P_0$-Donsker and
$\sup_{h\in\mathcal{H}}P_0|f_j-f_{0,j}|^2\rightarrow0$, as
$|\theta-\theta_0|+\|\eta-\eta_0\|\rightarrow0$,
for every $f_j\in\mathcal{F}_j(\delta_0),j=1,2$, where
$f_{0,1}=\dot{\ell}_{\theta_0,\eta_0}$ and $f_{0,2}=B_0h-P_0B_0h$,
(2) $\mathcal{F}_k(\delta_0),k=1,2$, have integrable envelopes.
\end{cond}

%
\begin{cond}
\label{condwlereg3}
The map
$\Psi= (\Psi_1,\Psi_2)\dvtx\Theta\times H\mapsto\mathbb
{R}^p\times
\ell^{\infty}(\mathcal{H})$
with components
\begin{eqnarray*}
\Psi_1(\theta,\eta) &\equiv& P_0\Psi_{N,1}(
\theta,\eta) =P_0\dot{\ell}_{\theta,\eta},
\\
\Psi_2(\theta,\eta)h &\equiv& P_0\Psi_{N,2}(
\theta,\eta)
=P_0B_{\theta,\eta}h-P_{\theta,\eta}B_{\theta,\eta}h,\qquad
h\in\mathcal{H},
\end{eqnarray*}
has a continuously invertible Fr$\acute{\mbox{e}}$chet derivative map
$\dot{\Psi}_0=(\dot{\Psi}_{11},\dot{\Psi}_{12},\dot{\Psi
}_{21},\dot{\Psi}_{22})$ at $(\theta_0,\eta_0)$ given by
$\dot{\Psi}_{ij}(\theta_0,\eta_0)h=P_0(\dot{\psi}_{i,j,\theta
_0,\eta_0,h})$, $i,j\in\{1,2\}$ in terms of $L_2(P_0)$ derivatives of
$\psi_{1,\theta,\eta,h}=\dot{\ell}_{\theta,\eta}$ and
$\psi_{2,\theta,\eta,h}=B_{\theta,\eta}h-P_{\theta,\eta
}B_{\theta,\eta}h$; that is,
\begin{eqnarray*}
\sup_{h\in\mathcal{H}} \bigl[P_0\bigl\{\psi_{i,\theta,\eta_0,h}-
\psi_{i,\theta_0,\eta_0,h} -\dot{\psi}_{i1,\theta_0,\eta
_0,h}(\theta
-\theta_0)\bigr
\}^2\bigr]^{1/2} &=& o\bigl(\|\theta-\theta_0
\|\bigr),
\\
\sup_{h\in\mathcal{H}} \bigl[P_0\bigl\{\psi_{i,\theta_0,\eta,h}-
\psi_{i,\theta_0,\eta_0,h} -\dot{\psi}_{i2,\theta_0,\eta
_0,h}(\eta-\eta
_0)\bigr
\}^2\bigr]^{1/2} &=& o\bigl(\|\eta-\eta_0\|\bigr).
\end{eqnarray*}
Furthermore, $\dot{\Psi}_0$ admits a partition
\[
(\theta-\theta_0,\eta-\eta) \mapsto%
\pmatrix{ \dot{
\Psi}_{11} & \dot{\Psi}_{12}
\cr
\dot{\Psi}_{21} &
\dot{\Psi}_{22}
\cr
} \pmatrix{ \theta-\theta_0
\cr
\eta-
\eta_0
\cr
},
\]
where
\begin{eqnarray*}
\dot{\Psi}_{11}(\theta-\theta_0) &=& -P_{\theta_0,\eta_0}
\dot{\ell}_{\theta_0,\eta_0}\dot{\ell}_{\theta
_0,\eta_0}^T(\theta-
\theta_0),
\\
\dot{\Psi}_{12}(\eta-\eta_0)&=&-\int B^*_{\theta_0,\eta_0}
\dot{\ell}_{\theta_0,\eta_0}\,d(\eta-\eta_0),
\\
\dot{\Psi}_{21}(\theta-\theta_0)h &=&-P_{\theta_0,\eta
_0}B_{\theta
_0,\eta_0}h
\dot{\ell}_{\theta
_0,\eta
_0}^T(\theta-\theta_0),
\\
\dot{\Psi}_{22}(\eta-\eta_0)h&=&-\int
B^*_{\theta_0,\eta
_0}B_{\theta
_0,\eta_0}h d(\eta-\eta_0)
\end{eqnarray*}
and $B^*_{\theta_0,\eta_0}B_{\theta_0,\eta_0}$ is continuously invertible.
\end{cond}

Let $\tilde{I}_0=P_0[(I-B_0(B^*_0B_0)^{-1}B_0^*)\dot{\ell}_0\dot
{\ell}_0^T]$
be the efficient information for $\theta$ and
$\tilde{\ell}_0=\tilde{I}_0^{-1}(I-B_0(B^*_0B_0)^{-1}B_0^*)\dot
{\ell}_0$
be the efficient influence function for $\theta$ for the semiparametric
model with complete data.
%
%
\begin{theorem}
\label{thmzthm1}
Under Conditions~\ref{condest}--\ref{condwlereg3},
\begin{eqnarray*}
\sqrt{N}(\hat{\theta}_N-\theta_0) &=&
\sqrt{N}\mathbb{P}_N^\pi\tilde{\ell}_0+o_{P^*}(1)
\rightsquigarrow Z \sim N_p(0,\Sigma),
\\
\sqrt{N}(\hat{\theta}_{N,\#}-\theta_0) &=& \sqrt{N}
\mathbb{P}_N^{\pi,\#}\tilde{\ell}_0+o_{P^*}(1)
\rightsquigarrow Z_{\#} \sim N_p(0,\Sigma_{\#}),
\end{eqnarray*}
where $\#\in\{e,c,\mathrm{mc},\mathrm{cc}\}$,
%
%
\begin{eqnarray}
\label{WLEasympVar}
\Sigma&\equiv& I_0^{-1} + \sum
_{j=1}^J\nu_j\frac{1-p_j}{p_j}
\operatorname{Var}_{0|j}(\tilde{\ell}_0),
\\
\label{WLECalasympVar}
\Sigma_{\#} &\equiv& I_0^{-1} + \sum
_{j=1}^J\nu_j\frac
{1-p_j}{p_j}
\operatorname{Var}_{0|j}\bigl((I-Q_{\#})\tilde{\ell}_0\bigr)
%
\end{eqnarray}
and (recall Conditions~\ref{condest}
and~\ref{condcal})
\begin{eqnarray*}
Q_ef&\equiv&P_0\bigl[\pi^{-1}_0(V)f
\dot{G}_e\bigl(Z^T\alpha_0
\bigr)Z^T\bigr]S^{-1}_0\bigl(1-
\pi_0(V)\bigr)^{-1}\dot{G}_e
\bigl(Z^T\alpha_0\bigr)Z,
\\
Q_{c}f&\equiv&P_0\bigl[f Z^T\bigr]\bigl
\{P_0Z^{\otimes2}\bigr\}^{-1}Z,
\\
Q_{\mathrm{mc}}f&\equiv&P_0\bigl[\bigl(\pi_0^{-1}(V)-1
\bigr)f Z^T\bigr]\bigl\{P_0\bigl[\bigl(
\pi_0^{-1}(V)-1\bigr)Z^{\otimes2}\bigr]\bigr
\}^{-1}Z,
\\
Q_{\mathrm{cc}}f&\equiv&P_0\bigl[\bigl(\pi_0^{-1}(V)-1
\bigr)f (Z-\mu_Z)^T\bigr]\bigl\{P_0\bigl[
\bigl(\pi_0^{-1}(V)-1\bigr) (Z-\mu_Z)^{\otimes2}
\bigr]\bigr\}^{-1}\\
&&{}\times(Z-\mu_Z).
\end{eqnarray*}
\end{theorem}

%
\begin{remark}
Our conditions in Theorem~\ref{thmzthm1} are the same as those in
\cite{MR2391566} except the integrability condition.
Our Condition~\ref{condwlereg2}(2) requires existence of integrable
envelopes for class of scores while the condition $(\mbox{A}1^*)$ in
\cite{MR2391566} requires square integrable envelopes.
Note that this integrability condition is required only for the WLE
with adjusted weights, as in~\cite{MR2325244}.
\end{remark}

%
\begin{remark}
As can be seen from the definition of $Q_\#$, the choice of $G$ in calibration
does not affect the asymptotic variances while $G_e$ in the method of
estimated weights does
affect the asymptotic variance.
\end{remark}

\subsection{Nonregular rate for a nuisance parameter}\label{sec3.3}
For\vspace*{1pt} $\underline{h}=(h_1,\ldots,h_p)^T$ with $h_k\in H$,
$k=1,\ldots,p$, let
$B_{\theta,\eta}[\underline{h}]=(B_{\theta,\eta}h_1,\ldots,B_{\theta,\eta}h_p)^T$.
We assume the following conditions.

%
\begin{cond}[(Consistency and rate of convergence)]
\label{condwlenonreg1}
An estimator $(\hat{\theta}_N,\hat{\eta}_N)$ of $(\theta_0,\eta
_0)$ satisfies
$|\hat{\theta}_N-\theta_0|= o_P(1)$, and
$\|\hat{\eta}_N-\eta_0\|= O_P(N^{-\beta})$ for some
$\beta>0$.
\end{cond}

%
\begin{cond}[(Positive information)]
\label{condwlenonreg2}
There is an $\underline{h}^*= (h_1^*,\ldots,h_p^*)$, where
$h_k^*\in H$ for $ k=1,\ldots, p$, such that
\[
P_0 \bigl\{ \bigl( \dot{\ell}_{0} -B_{0}
\bigl[\underline{h}^*\bigr] \bigr) B_{0}h \bigr\}= 0 \qquad\mbox
{for all } h
\in H.
\]
The efficient information
$I_0\equiv P_0(\dot{\ell}_{0}-B_{0}[\underline{h}^*])^{\otimes2} $
for $\theta$ for the semiparametric model with complete data
is finite and nonsingular.
Denote the efficient influence function for the semiparametric model
with complete data by $\tilde{\ell}_0 \equiv I_0^{-1}(\dot{\ell
}_{0}-B_{0}[\underline{h}^*])$.
\end{cond}
%
%
\begin{cond}[(Asymptotic equicontinuity)]
\label{condwlenonreg3}
(1) For any $\delta_N\downarrow0$ and $C>0$,
\begin{eqnarray*}
\sup_{|\theta-\theta_0|\leq\delta_N,\|\eta-\eta
_0\|
\leq CN^{-\beta}} \bigl|\mathbb{G}_N(\dot{
\ell}_{\theta,\eta}-\dot{\ell}_{0})\bigr| &=& o_P(1),
\\
\sup_{|\theta-\theta_0|\leq\delta_N,\|\eta-\eta
_0\|
\leq CN^{-\beta}} \bigl|\mathbb{G}_N(B_{\theta,\eta}-B_{0})
\bigl[\underline{h}^*\bigr]\bigr| &=& o_P(1).
\end{eqnarray*}

(2) There exists a $\delta>0$ such that the classes
$\{\dot{\ell}_{\theta,\eta}\dvtx|\theta-\theta_0| +\|\eta
-\eta_0\|
\leq\delta\}$
and
$\{B_{\theta,\eta} [\underline{h}^* ]\dvtx|\theta-\theta_0|+\|
\eta-\eta_0\|\leq\delta\}$
are $P_0$-Glivenko--Cantelli and have integrable envelopes.
Moreover, $\dot{\ell}_{\theta,\eta}$ and $B_{\theta,\eta
}[\underline
{h}^*]$ are
continuous with respect to $(\theta,\eta)$ either pointwise or in $L_1(P_0)$.
\end{cond}
%
%
\begin{cond}[(Smoothness of the model)]
\label{condwlenonreg4}
For some $\alpha>1$ satisfying $\alpha\beta>1/2$ and for $(\theta,\eta)$
in the neighborhood $\{(\theta,\eta)\dvtx|\theta-\theta_0|\leq
\delta
_N,\|\eta-\eta_0\|\leq CN^{-\beta}\}$,
\begin{eqnarray*}
&&\bigl|P_0\bigl\{\dot{\ell}_{\theta,\eta}-\dot{\ell}_{0} +
\dot{\ell}_{0}\bigl(\dot{\ell}_{0}^T(\theta-
\theta_0) +B_{0}[\eta-\eta_0]\bigr)\bigr\}\bigr|
\\
&&\qquad=o\bigl(|\theta-\theta_0|\bigr)+O\bigl(\|\eta-\eta_0
\|^{\alpha
}\bigr),
\\
&&\bigl|P_0\bigl\{(B_{\theta,\eta}-B_{0})\bigl[
\underline{h}^*\bigr] +B_{0}\bigl[\underline{h}^*\bigr]\bigl(\dot{
\ell}_{0}^T (\theta-\theta_0)+B_{0}[
\eta-\eta_0]\bigr)\bigr\}\bigr|
\\
&&\qquad= o\bigl(|\theta-\theta_0|\bigr)+O\bigl(\|\eta-\eta_0
\|^{\alpha}\bigr).
\end{eqnarray*}
\end{cond}

In the previous section, we required that the WLE
solves the weighted likelihood equations (\ref{eqnwlik}) for all
$h\in
\mathcal{H}$.
Here, we only assume that the WLE $(\hat{\theta}_N,\hat{\eta}_N)$
satisfies the weighted likelihood equations
%
%
\begin{eqnarray}
\label{eqnwlik2} \Psi_{N,1}^\pi(\theta,\eta,\alpha)&=&
\mathbb{P}_N^{\pi}\dot{\ell}_{\theta,\eta}
=o_{P^*}\bigl(N^{-1/2}\bigr),
\nonumber\\[-8pt]\\[-8pt]
\Psi_{N,2}^\pi(\theta,\eta,\alpha)\bigl[\underline{h}^*
\bigr] &=&\mathbb{P}_N^{\pi} B_{\theta,\eta}\bigl[
\underline{h}^*\bigr]=o_{P^*}\bigl(N^{-1/2}\bigr).\nonumber
\end{eqnarray}
The corresponding WLEs with adjusted weights, $(\hat{\theta}_{N,\#
},\hat{\eta}_{N,\#})$ with $\#\in\{e,c,\mathrm{mc},\mathrm{cc}\}$
satisfy (\ref{eqnwlik2})
with $\mathbb{P}_N^\pi$ replaced by $\mathbb{P}_N^{\pi,\#}$.

%
\begin{theorem}
\label{thmzthm2}
Suppose that the WLE is a solution of
(\ref{eqnwlik2}) where $\mathbb{P}_N^\pi$ may be replaced by
$\mathbb
{P}_N^{\pi,\#}$ with $\#\in\{e,c,\mathrm{mc},\mathrm{cc}\}$
for the estimators with adjusted weights.
Under Conditions~\ref{condest},~\ref{condcal} and \ref
{condwlenonreg1}--\ref{condwlenonreg4},
\begin{eqnarray*}
\sqrt{N}(\hat{\theta}_N-\theta_0) &=& \sqrt{N}
\mathbb{P}_N^\pi\tilde{\ell}_0+o_{P^*}(1)
\rightsquigarrow Z \sim N_p(0,\Sigma),
\\
\sqrt{N}(\hat{\theta}_{N,\#}-\theta_0) &=& \sqrt{N}
\mathbb{P}_N^{\pi,\#}\tilde{\ell}_0+o_{P^*}(1)
\rightsquigarrow Z_{\#} \sim N_p(0,\Sigma_{\#}),
\end{eqnarray*}
where $\Sigma$ and $\Sigma_{\#}$ are as
defined in (\ref{WLEasympVar}) and (\ref{WLECalasympVar}) of Theorem
\ref{thmzthm1},
but now $I_0$ and $\tilde{\ell}_0$ are defined in Condition
\ref{condwlenonreg2}, and $Q_{\#}$ are defined in Theorem~\ref{thmzthm1}.
\end{theorem}

%
\begin{remark}
Our conditions are identical to those of the $Z$-theorem of \cite
{MR1394975} 
except Condition~\ref{condwlenonreg3}(2).
This additional condition is not stringent for the following reasons.
First, the Glivenko--Cantelli condition is usually assumed to prove
consistency of estimators before deriving asymptotic distributions.
Second, a stronger $L_2(P_0)$-continuity condition is standard
as is seen in Condition~\ref{condwlereg2} (see also Section 25.8 of
\cite{MR1652247}).
Note that the $L_1(P_0)$-continuity condition is only required for the
WLEs with adjusted weights.
\end{remark}

\subsection{Comparisons of methods}\label{sec3.4}
We compare asymptotic variances of five WLEs in view of improvement
by adjusting weights and change of designs.
We also include in comparison special cases of adjusting weights
involving stratum-wise adjustment.

\subsubsection{Stratified Bernoulli sampling}\label{sec3.4.1}
We first give a statement of the result corresponding to Theorem
\ref{thmzthm1} for stratified Bernoulli sampling where all
subjects are\vadjust{\goodbreak} independent with the sampling probability $p_j$ if
$V\in\mathcal{V}_j$ and $\hat{\theta}_N^{\mathrm{Bern}}$ and
$\hat{\theta}_{N,\#}^{\mathrm{Bern}}$ with
$\#\in\{e,c,\mathrm{mc},\mathrm{cc}\}$ are the corresponding WLE and
WLEs with adjusted weights.
%
%
\begin{theorem}
\label{thmzthmBern}
Suppose Conditions~\ref{condest} [except Condition~\ref{condest}\textup{(f)}] and
\ref{condcal} hold.
Let $\xi_i\in\{0,1\}$ and $\xi$ be i.i.d. with \mbox{$E[\xi|V]=\pi
_0(V)=\sum
_{j=1}^Jp_jI(V\in\mathcal{V}_j)$}.

(1) Suppose that the WLE is a solution of
(\ref{eqnwlik2}) where $\mathbb{P}_N^\pi$ may be replaced by
$\mathbb
{P}_N^{\pi,\#}$ with $\#\in\{e,c,\mathrm{mc},\mathrm{cc}\}$
for the estimators with adjusted weights.
Under the same conditions as in Theorem~\ref{thmzthm1},
\begin{eqnarray*}
\sqrt{N}\bigl(\hat{\theta}_N^{\mathrm{Bern}}-\theta_0
\bigr) &=& \sqrt{N}\mathbb{P}_N^\pi\tilde{
\ell}_0+o_{P^*}(1) \rightsquigarrow Z^{\mathrm{Bern}} \sim
N_p\bigl(0,\Sigma^{\mathrm{Bern}}\bigr),
\\
\sqrt{N}\bigl(\hat{\theta}_{N,\#}^{\mathrm{Bern}}-\theta_0
\bigr) &=& \sqrt{N}\mathbb{P}_N^{\pi,\#}\tilde{
\ell}_0+o_{P^*}(1) \rightsquigarrow Z_{\#}^{\mathrm{Bern}}
\sim N_p\bigl(0,\Sigma_{\#}^{\mathrm{Bern}}\bigr),
\end{eqnarray*}
where
%
%
\begin{eqnarray}
\label{WLEasympVarBer}
\Sigma^{\mathrm{Bern}} &\equiv& I_0^{-1} + \sum
_{j=1}^J\nu_j\frac
{1-p_j}{p_j}P_{0|j}(
\tilde{\ell}_0)^{\otimes2},
\\
\label{WLECalasympVarBer}
\Sigma_{\#}^{\mathrm{Bern}} &\equiv& I_0^{-1} +
\sum_{j=1}^J\nu_j
\frac
{1-p_j}{p_j}P_{0|j}\bigl((I-Q_{\#})\tilde{
\ell}_0\bigr)^{\otimes2}, %
%
\end{eqnarray}
where $Q_{\#}$ with $\#\in\{e,c,\mathrm{mc},\mathrm{cc}\}$ are defined in Theorem
\ref{thmzthm1}.

(2) Under the same conditions as in Theorem~\ref{thmzthm2},
the same conclusions in (1) hold with $I_0$ and $\tilde{\ell}_0$
replaced by those defined in Condition~\ref{condwlenonreg2}.
\end{theorem}

Comparing the variance--covariance matrices in
Theorem~\ref{thmzthmBern} to those in
Theorems~\ref{thmzthm1} and~\ref{thmzthm2}, we obtain the following
corollary comparing designs.
All estimators have smaller variances under sampling without replacement.

%
\begin{cor}
\label{cordes}
Under the same conditions as in Theorem~\ref{thmzthmBern},
\begin{eqnarray*}
\Sigma&=&\Sigma^{\mathrm{Bern}} -\sum_{j=1}^J
\nu_j\frac{1-p_j}{p_j}\{ P_{0|j}\tilde{\ell}_0
\}^{\otimes2},
\\
\Sigma_{\#}&=&\Sigma^{\mathrm{Bern}}_{\#} -\sum
_{j=1}^J\nu_j\frac
{1-p_j}{p_j}\bigl\{
P_{0|j}(I-Q_{\#})\tilde{\ell}_0\bigr
\}^{\otimes2},\qquad \#\in\{e,c,\mathrm{mc},\mathrm{cc}\}. 
\end{eqnarray*}
%
\end{cor}

Variance formulas (\ref{WLECalasympVarBer}) with $\#\in\{e,\mathrm{mc},\mathrm{cc}\}$
except for the ordinary calibration have the following alternative
representations
which show the efficiency gains over the plain WLE under Bernoulli sampling.
%
%
\begin{cor}
\label{corberimp}
Under the same conditions as in Theorem~\ref{thmzthmBern},
\[
\Sigma_{\#}^{\mathrm{Bern}} = \Sigma^{\mathrm{Bern}} - \operatorname{Var}
\biggl(\frac{\xi-\pi_0(V)}{\pi_0(V)}Q_{\#
}\tilde{\ell}_0 \biggr),\qquad \#
\in\{e,\mathrm{mc},\mathrm{cc}\}. %
\]
%
\end{cor}
%

\subsubsection{Within-stratum adjustment of weights}\label{sec3.4.2}
Adjusting weights can be carried out in every stratum.
This is proposed by Breslow et al.~\cite{Breslow2009a,Breslow2009b}
for ordinary calibration.
Consider\vspace*{1pt} calibration on $\tilde{Z}$ where $\tilde{Z}\equiv
(Z^{(1)},\ldots,Z^{(J)})^T$ with $Z^{(j)}\equiv I(V\in\mathcal{V}_j)Z^T$.
The calibration equation (\ref{eqncaleqn2}) becomes
\[
\frac{1}{N}\sum_{i=1}^N
\frac{\xi_iG_c(\tilde{Z}_i;\alpha)}{\pi_0(V_i)}Z_iI(V_i\in
\mathcal{V}_j)
=\frac{1}{N}\sum_{i=1}^NZ_iI(V_i
\in\mathcal{V}_j),\qquad j=1,\ldots,J,
\]
where $\alpha\in\mathbb{R}^{Jk}$.
We call this special case \textit{within-stratum calibration}.
We define \textit{within-stratum modified and centered calibration}
analogously.

We also call estimated weights carried out within stratum
\textit{within-stratum estimated weights}. Recall that $Z$ in
estimated weights
contains the membership indicators for the strata and the rest are
other auxiliary variables, say $Z^{[2]}$. Within-stratum estimated
weights uses $\tilde{Z}\equiv(Z^{(1)},\ldots,Z^{(J)})^T$ where
$Z^{(j)}\equiv I(V\in\mathcal{V}_j)(Z^{[2]})^T$ with 1 included in
$Z^{[2]}$. The ``true'' parameter $\tilde{\alpha}_0$ has zero for all
elements except having $G^{-1}_e(p_j)$ for the element corresponding to
$I(V\in\mathcal{V}_j)$, $j=1,\ldots,J$.

The following corollary summarizes within-stratum adjustment of weights under
stratified Bernoulli sampling and sampling without replacement.
All methods achieve improved efficiency over the plain WLE under
Bernoulli sampling while centered calibration is the only method to
yield a
guaranteed improvement under sampling without replacement.
This is because centering yields the $L_2^0(P_{0|j})$-projection
suitable for
the conditional variances in (\ref{WLECalasympVar}) while noncentering
results in the $L_2(P_{0|j})$-projection for the conditional
expectations in (\ref{WLECalasympVarBer}).
%
%
\begin{cor}
\label{corst}
(1) (Bernoulli) Under the same conditions as in Theorem~\ref{thmzthmBern}
with $Z$ replaced by $\tilde{Z}$ and $\alpha_0$ replaced by
$\tilde{\alpha}_0$ for within-stratum estimated weights,
%
%
\begin{equation}\label{eqnCalWImproveBern}
\Sigma_{\#}^{\mathrm{Bern}} = \Sigma^{\mathrm{Bern}} - \sum
_{j=1}^J\nu_j\frac{1-p_j}{p_j}P_{0|j}
\bigl(Q_{\#
}^{(j)}\tilde{\ell}_0
\bigr)^{\otimes2}, %
\end{equation}
where $\#\in\{e,c,\mathrm{mc},\mathrm{cc}\}$ and
\begin{eqnarray*}
Q_e^{(j)}f &\equiv&P_{0|j} \bigl[f
\dot{G}_e\bigl(\tilde{Z}^T\tilde{\alpha}_0
\bigr) \bigl(Z^{[2]}\bigr)^T \bigr] \bigl\{P_{0|j}
\dot{G}^2_e\bigl(\tilde{Z}^T\tilde{\alpha
}_0\bigr) \bigl(Z^{[2]}\bigr)^{\otimes2} \bigr
\}^{-1}
\\
&&{}\times\dot{G}_e\bigl(\tilde{Z}^T\tilde{
\alpha}_0\bigr)I(V\in\mathcal{V}_j)Z^{[2]},
\\
Q_{c}^{(j)}f &\equiv&P_{0|j}\bigl[f Z^T
\bigr]\bigl\{P_{0|j}\bigl[Z^{\otimes2}\bigr]\bigr\}^{-1}I(V
\in\mathcal{V}_j)Z,
\\
Q_{\mathrm{mc}}^{(j)}f &\equiv& Q_{c}^{(j)}f,
\\
Q_{\mathrm{cc}}^{(j)}f&\equiv&P_{0|j}\bigl[f (Z-
\mu_{Z,j})^T\bigr]\bigl\{P_{0|j}\bigl[(Z-
\mu_{Z,j})^{\otimes2}\bigr]\bigr\}^{-1}I(V\in
\mathcal{V}_j) (Z-\mu_{Z,j})
\end{eqnarray*}
with $\mu_{Z,j}\equiv E [I(V\in\mathcal{V}_j)Z]$ for $j=1,\ldots,J$.\vadjust{\goodbreak}

(2) (without replacement)
Under the same conditions as in Theorems~\ref{thmzthm1} or
\ref{thmzthm2} with $Z$ is replaced by $\tilde{Z}$,
%
%
\begin{equation}\label{eqnCentCalWImproveSWOR}
\Sigma_{\mathrm{cc}} = \Sigma- \sum_{j=1}^J
\nu_j\frac{1-p_j}{p_j}\operatorname{Var}_{0|j}
\bigl(Q_{\mathrm{cc}}^{(j)}\tilde{\ell}_0 \bigr).
\end{equation}
\end{cor}

\subsubsection{Comparisons}\label{sec3.4.3}
We summarize Corollaries~\ref{cordes}--\ref{corst}.
Every method of adjusting weights improves efficiency over the plain
WLE in a certain design and with a certain range of adjustment
of weights (within-stratum or ``across-strata'' adjustment).
However, particularly notable among all methods is centered calibration.
While other methods gain efficiency only under Bernoulli sampling,
centered calibration improves efficiency over the plain WLE under both
sampling schemes.
There is no known method of ``across-strata'' adjustment
that is guaranteed to gain efficiency over the plain WLE under
stratified sampling without replacement.

There are close connections among all methods.
When the auxiliary variables have mean zero, centered
and modified calibrations are essentially the same.
The ordinary and modified calibrations give the same
asymptotic variance when carried out stratum-wise.
For $Z$ and $\alpha_0$ defined for estimated weights,\vspace*{1pt} estimated weights
and modified calibration based on $(1-\pi_0(V))^{-1}\dot
{G}_e(Z^T\alpha_0)Z$ performs the same way.
Similarly\vspace*{2pt} within-stratum estimated weights with $\tilde{Z}$ and
$\tilde{\alpha}_0$ is as good as within-stratum calibration based on
$\dot{G}_e(\tilde{Z}^T\tilde{\alpha_0})\tilde{Z}$.

As seen in the relationship among methods, there is
no single method superior to others in each situation.
In fact, performance depends on choice and transformation
of auxiliary variables, the true distribution $P_0$ and the design.
For our ``without replacement'' sampling scheme, within-stratum centered
calibration is
the only method guaranteed to gain efficiency while other methods
may perform even worse than the plain WLE.

\section{Examples}\label{sec4}
For asymptotic normality of WLEs, consistency and rate of convergence
need to be established first to apply our $Z$-theorems in Section~\ref{sec3}.
To this end, general results on IPW empirical processes discussed in
the next section will be useful.
We illustrate this in the Cox models
with right censoring and interval censoring under two-phase sampling.

Let $T\sim F$ be a failure time, and $X$ be a vector of covariates
with bounded supports in the regression model.
The Cox proportional hazards model~\cite{MR0341758} specifies the relationship
\[
\Lambda(t|x)=\exp\bigl(\theta^Tx\bigr)\Lambda(t),
\]
where $\theta\in\Theta\subset\mathbb{R}^p$ is the regression parameter,
$\Lambda\in H$ is the (baseline) cumulative hazard function.
Here the space $H$ for the nuisance parameter $\Lambda$ is the set of
nonnegative,
nondecreasing cadlag functions defined on the positive line.
The true parameters are $\theta_0$ and $\Lambda_0$.\vadjust{\goodbreak}

In addition to $X$, let $U$ be a vector of auxiliary variables collected
at phase I which are correlated with the covariate $X$.
For simplicity of notation, we assume that the covariates $X$ are only
observed for the subject sampled at phase II.
Thus, if some of the coordinates of $X$ are available at phase I,
then we include an identical copy of those coordinates of $X$ in the
vector of $U$.

\subsection{Cox model with right censored data}\label{sec4.1}
Under right censoring, we only observe the minimum
of the failure time $T$ and the censoring time $C\sim G$.
Define the observed time $Y=T\wedge C$ and the censoring indicator
$\Delta=I(T\leq C)$.
The phase I data is $V=(Y,\Delta,U)$, and the observed data is
$(Y,\Delta,\xi X,U,\xi)$ where $\xi$ is the sampling indicator.
With right censored data and complete data, the theory for maximum likelihood
estimators in the Cox model has received several treatments; the one we
follow most
closely here is that of~\cite{MR1652247}. For the Cox model with
case-cohort data, see
\cite{MR924857} and for treatments with even more general designs
\cite{MR1158522} and~\cite{MR1766826}. Here, for both sampling without
replacement
and Bernoulli sampling, we continue the developments of \cite
{MR2325244,MR2391566}.
We assume the following conditions:

%
\begin{cond}
\label{condcoxr1}
The finite-dimensional parameter space $\Theta$ is compact
and contains the true parameter $\theta_0$ as an interior point.
\end{cond}
%
%
\begin{cond}
\label{condcoxr2}
The failure time $T$ and the censoring time $C$ are conditionally
independent given $X$, and that there is
$\tau>0$ such that $P(T>\tau)>0$ and $P(C\geq\tau)=P(C=\tau)>0$.
Both $T$ and $C$ have continuous conditional densities given the
covariates $X=x$.
\end{cond}
%
%
\begin{cond}
\label{condcoxr3}
The covariate $X$ has bounded support. For any measurable function
$h$, $P(X\neq h(Y))>0$.
\end{cond}

Let $\lambda(t)=(d/dt)\Lambda(t)$ be the baseline hazard function.
With complete data, the density of $(Y,\Delta,X)$ is
\[
p_{\theta,\Lambda}(y,\delta,x) = \bigl\{\lambda(y)e^{\theta^Tx
-\Lambda
(y)e^{\theta^Tx}}(1-G) (y|x)
\bigr\}^\delta\bigl\{e^{-\Lambda(y)e^{\theta^Tx}}g(y|x)\bigr\}
^{1-\delta}p_X(x),
\]
where $p_X$ is the marginal density of $X$ and $g(\cdot|x)$ is the
conditional density of $C$ given $X=x$.
The score for $\theta$ is given by $\dot{\ell}_{\theta,\Lambda} (y,
\delta, x)= x \{ \delta- e^{\theta^Tx}\Lambda(y) \}$,
and the score operator $B_{\theta,\Lambda}\dvtx\mathcal{H}\mapsto
L_2(P_{\theta,\Lambda})$
is defined on the unit ball $\mathcal{H}$ in the space $BV[0,\tau]$
such that
$B_{\theta,\Lambda}h (y, \delta, x) = \delta h(y) -e^{\theta
^Tx}\int_{[0,y]}h\,d\Lambda$.
Because the likelihood based on the density above does not yield the
MLE for complete data, we define the log likelihood for
one observation for complete data by $\ell_{\theta,\Lambda} (y,
\delta, x)
= \log\{(e^{\theta^Tx}\Lambda\{y\})^\delta e^{-\Lambda(y)e^{\theta
^Tx}}\}$
where $\Lambda\{t\}$ is the (point) mass of $\Lambda$ at $t$.
Then maximizing the weighted log likelihood $\mathbb{P}_N^\pi\ell
_{\theta,\Lambda}$
reduces to solving the system of equations $\mathbb{P}_N^\pi\dot
{\ell
}_{\theta,\Lambda}=0$
and $\mathbb{P}_N^\pi B_{\theta,\Lambda}h=0$ for every $h\in
\mathcal{H}$.
The efficient score for $\theta$ for complete data is given by
\[
\ell^*_{\theta_0,\Lambda_0} (y, \delta, x) =\delta\bigl(x-(M_1/M_0)
(y)\bigr) -e^{\theta_0^Tx}\int_{[0,y]}\delta
\bigl(x-(M_1/M_0) (t) \bigr)\,d\Lambda_0(t),
\]
and the efficient information for $\theta$ for complete data is
\[
\tilde{I}_{\theta_0,\Lambda_0} = E \bigl[ \bigl(\ell^*_{\theta
_0,\Lambda_0}
\bigr)^{\otimes
2} \bigr] =Ee^{\theta_0^TX}\int_0^\tau
\biggl(X-\frac{M_1}{M_0}(y) \biggr)^{\otimes
2}(1-G) (y|X)\,d
\Lambda_0(y),
\]
where $M_k(s) = P_{\theta_0,\Lambda_0}[X^ke^{\theta_0^TX}I(Y\geq s)]$,
$k= 0, 1$.

%
\begin{theorem}[(Consistency)]
\label{thmwlecons} Under Conditions~\ref{condest},~\ref{condcal},
\ref{condcoxr1}--\ref{condcoxr3}, the WLEs are consistent for
$(\theta_0,\Lambda_0)$.
\end{theorem}

\begin{pf}
This proof follows along the lines of the proof given by \cite
{MR1915446}, but with
the usual empirical measure replaced by the IPW empirical measure (with
adjusted weights),
and by use of Theorem~\ref{thmfpsPiGC}.
For details see~\cite{Saegusa-Wellner2012supp}.
\end{pf}

Our $Z$-theorem (Theorem~\ref{thmzthm1}) yields asymptotic normality
of the WLEs.

%
\begin{theorem}[(Asymptotic normality)]
\label{thmanright}
Under Conditions~\ref{condest},~\ref{condcal},\break
\mbox{\ref{condcoxr1}--\ref{condcoxr3}},
\begin{eqnarray*}
\sqrt{N}(\hat{\theta}_N-\theta_0) &=&\sqrt{N}
\mathbb{P}_N^\pi\tilde{\ell}_{\theta_0,\Lambda_0}+o_{P^*}(1)
\rightarrow_d N (0,\Sigma),
\\
\sqrt{N}(\hat{\theta}_{N,\#}-\theta_0) &=&\sqrt{N}
\mathbb{P}_N^{\pi,\#}\tilde{\ell}_{\theta_0,\Lambda
_0}+o_{P^*}(1)
\rightarrow_d N (0,\Sigma_{\#} ), 
\end{eqnarray*}
where $\#\in\{e,c,\mathrm{mc},\mathrm{cc}\}$,
$\tilde{\ell}_{\theta_0,\Lambda_0}=I^{-1}_{\theta_0,\Lambda
_0}\ell^*_{\theta_0,\Lambda_0}$
is the efficient influence function for $\theta$ for complete data,
and $\Sigma$ and $\Sigma_{\#}$
are given in Theorem~\ref{thmzthm1}.
\end{theorem}

\begin{pf}
We verify the conditions of Theorem~\ref{thmzthm1}.
Condition~\ref{condwlereg1} holds by Theorem~\ref{thmwlecons}.
Conditions~\ref{condwlereg2} and~\ref{condwlereg3} hold under the
present hypotheses as was shown in
\cite{MR1652247}, Section 25.12. 
\end{pf}
For variance estimation regarding
$\hat{\theta}_N$,
$\hat{I}_N\equiv\mathbb{P}_N^\pi(\ell^*_{\hat{\theta}_N,\hat
{\Lambda
}_N})^{\otimes2}$
can be used to estimate $I_0$.
Letting
$\hat{\tilde{\ell}}_0\equiv\hat{I}_N^{-1} \ell^*_{\hat{\theta
}_N,\hat
{\Lambda}_N}$,
we can estimate $\operatorname{Var}_{0|j}\tilde{\ell}_0$ by
$\hat{P}_{j}\tilde{\ell}_0^{\otimes2}-\{\hat{P}_j\tilde{\ell}_0\}
^{\otimes2}$
where $\hat{P}_j\tilde{\ell}_0\equiv\mathbb{P}_N^\pi\hat{\tilde
{\ell
}}_0I(V\in\mathcal{V}_j)$
and
$\hat{P}_j\tilde{\ell}_0^{\otimes2}
\equiv\mathbb{P}_N^\pi\hat{\tilde{\ell}}{}^{\otimes2}_0 I(V\in
\mathcal{V}_j)$.
The other four cases are similar.

\subsection{Cox model with interval censored data}\label{sec4.2}
Let $Y$ be a censoring time that is assumed to be conditionally
independent of a failure time $T$ given a covariate vector $X$.
Under the case 1 interval censoring, we do not observe
$T$ but $(Y,\Delta)$ where $\Delta\equiv I(T\leq Y)$.
The phase I data is $V=(Y,\Delta,U)$ and the observed data is
$(Y,\Delta,\xi X,U,\xi)$ where $\xi$ is the sampling indicator.
In the case of complete data, maximum likelihood estimates for this
model were studied
by Huang~\cite{MR1394975}.
For a generalized version of this model and
two-phase data with Bernoulli sampling, weighted likelihood estimates
with and without estimated
weights have recently been studied by Li and Nan~\cite{MR2860827}.
Here we treat two-phase data under sampling without replacement at
phase II and with
both estimated weights and calibration.

With complete data, the log likelihood for one observation is given by
\begin{eqnarray*}
\ell(\theta,F) &\equiv&\delta\log\bigl\{1-\overline{F}(y)^{\exp
(\theta
^Tx)}\bigr\}
+(1-\delta)\log\overline{F}(y)^{\exp(\theta^Tx)}
\\
&\equiv&\delta\log\bigl\{1-e^{-\Lambda(y)\exp(\theta^Tx)}
\bigr\}-(1-\delta)e^{\theta^Tx }\Lambda(y) \equiv\ell(\theta,\Lambda),
\end{eqnarray*}
where $\overline{F}\equiv1-F = e^{-\Lambda}$.
The score for $\theta$ and the score operator $B_{\theta,\Lambda}$ for
$\Lambda$ for complete data are
$\dot{\ell}_{\theta,\Lambda}=x\exp(\theta^Tx)\Lambda(y)(\delta
r(y,x;\theta,\Lambda)-(1-\delta))$
and
$B_{\theta,\Lambda}[h]
=\exp(\theta^Tx)h(y) \{\delta r(y,x;\theta,\Lambda)-(1-\delta
)
\}$
where
$r(y,x;\theta,\Lambda) =\break\exp(-e^{\theta^Tx}\Lambda(y))/\{1-\exp
(-e^{\theta^Tx}\Lambda(y))\}$.
The efficient score for $\theta$ for complete data is given by
\[
\ell^*_{\theta_0,\Lambda_0} =e^{\theta_0^Tx}Q(y,\delta,x;\theta_0,
\Lambda_0)\Lambda_0(y) \biggl\{x-\frac{E[Xe^{2\theta
_0^TX}O(Y|X)|Y=y]}{E[e^{2\theta
_0^TX}O(Y|X)|Y=y]} \biggr
\},
\]
where $Q(y,\delta,x;\theta,\Lambda) =\delta r(y,x;\theta,\Lambda
)-(1-\delta)$ and
$O(y|x)= r(y,x; \theta_0, \Lambda_0) $.
The efficient information for $\theta$ for complete data
$\tilde{I}_{\theta_0,\Lambda_0} = E[(\ell^*_{\theta_0,\Lambda
_0})^{\otimes2}]$
is given by $\tilde{I}_{\theta_0,\Lambda_0}
=E[R(Y,X)\{X-E[XR(Y,X)|Y]/E[R(Y,X)|Y]\}]$
where\break $R(Y,X)=\Lambda_0^2(Y|X)O(Y|X)$.
See~\cite{MR1394975} 
for further details.

We impose the same assumptions made for complete data in \cite
{MR1394975}. 
%
%
\begin{cond}
\label{condcoxi1}
The finite-dimensional parameter space $\Theta$ is compact
and contains the true parameter $\theta_0$ as its interior point.
\end{cond}
%
%
\begin{cond}
\label{condcoxi2}
(a) The covariate $X$ has bounded support; that is, there exists $x_0$
such that $|X|\leq x_0$ with probability 1.
(b) For any $\theta\neq\theta_0$, the probability $P(\theta^TX\neq
\theta_0^T X)>0$.
\end{cond}
%
%
\begin{cond}
\label{condcoxi3}
$F_0(0)=0$. Let $\tau_{F_0}=\inf\{t\dvtx F_0(t)=1\}$.
The support of $Y$ is an interval $S[Y]=[l_Y,u_Y]$ and $0< l_Y\leq
u_Y< \tau_{F_0}$.
\end{cond}
%
%
\begin{cond}
\label{condcoxi4}
The cumulative hazard function $\Lambda_0$ has strictly positive
derivative on
$S[Y]$, and the joint function $G(y,x)$ of $(Y,X)$ has bounded second
order (partial)
derivative with respect to $y$.
\end{cond}

\subsubsection{Consistency}\label{sec4.2.1}
The characterization of WLEs $(\hat{\theta}_N,\hat{\Lambda}_N)$
and $(\hat{\theta}_{N,\#},\break\hat{\Lambda}_{N,\#})$ with
$\#\in\{e,c,\mathrm{mc},\mathrm{cc}\}$ maximizing $\mathbb{P}_N^\pi\ell(\theta,\Lambda)$ or
$\mathbb{P}_N^{\pi,\#} \ell(\theta,\Lambda)$ is given in~\cite
{Saegusa-Wellner2012supp}, Lemma A.5.
We prove consistency of the WLEs in the metric given by
$d((\theta_1,\Lambda_1),(\theta_2,\Lambda_2))
\equiv\|\theta_1-\theta_2\|+ \|\Lambda_1-\Lambda_2
\|_{P_Y}$,
where \mbox{$\|\cdot\|$} is the Euclidean metric and
$ \|\Lambda_1-\Lambda_2 \|_{P_Y}^2
=\int(\Lambda_1(y)-\Lambda_2(y) )^2\,dP_Y,
$
and $P_Y$ is the marginal probability measure of the censoring variable $Y$.

%
\begin{theorem}[(Consistency)]
\label{thmconscurrent} Under Conditions~\ref{condest},~\ref{condcal},
\ref{condcoxi1}--\ref{condcoxi4}, the WLEs are consistent in the
metric $d$.
\end{theorem}

\begin{pf}
We only prove consistency for the WLE.
Proofs for the other four estimators are similar.

Let $\tilde{H}$ be the set of all subdistribution functions defined on
$[0,\infty]$.
We denote the WLE of $F$ as $\hat{F}_N=1-e^{-\hat{\Lambda}_N}$.
Define the set $\mathcal{F}$ of functions by
\[
\mathcal{F} \equiv\bigl\{f(\theta,F)=\delta\bigl(1-\overline
{F}(y)^{\exp(\theta^Tx)}
\bigr) +(1-\delta) \overline{F}(y)^{\exp(\theta^Tx)}\dvtx\theta
\in\Theta, F
\in\tilde{H}
\bigr\}.
\]
Boundedness of $X$ and compactness of $\Theta\subset\mathbb{R}^p$
imply that the set $\{e^{\theta^Tx}\dvtx\theta\in\Theta\}$ is
Glivenko--Cantelli.
The set $\tilde{H}$ is also Glivenko--Cantelli since it is a subset of the
set of bounded monotone functions.
Thus, it follows from boundedness of functions in $\mathcal{F}$
and the Glivenko--Cantelli preservation theorem~\cite{MR1857319} that
$\mathcal{F}$ is Glivenko--Cantelli.

Let $0<\alpha<1$ be a fixed constant.
It follows by concavity of the function $u\mapsto\log u$
and Jensen's inequality that
\begin{eqnarray*}
&&
P_{0}\bigl[\log
\bigl\{1+\alpha\bigl(f(\theta,F)/f(\theta_0,F_0)-1\bigr)
\bigr\}\bigr]\\
&&\qquad\leq\log\bigl(P_{0}\bigl[1+\alpha\bigl(f(\theta,F)/f(
\theta_0,F_0)-1\bigr)\bigr]\bigr)
\\
&&\qquad= \log\bigl(1-\alpha+\alpha P_{0}\bigl[f(\theta,F)/f(
\theta_0,F_0)\bigr] \bigr) \leq0,
\end{eqnarray*}
where the first equality holds if and only if $1+\alpha(f(\theta,F)/f(\theta_0,F_0)-1)$
is constant on $S[Y]$, in other words, $(\theta,F)=(\theta_0,F_0)$ on $S[Y]$
by the identifiability Condition~\ref{condcoxi2}.
Note also that by monotonicity of the logarithm
\begin{eqnarray*}
P_{0}\bigl[\log\bigl
\{1+\alpha\bigl(f(\theta,F)/f(\theta_0,F_0)-1\bigr)\bigr\}
\bigr] &\geq& P_{0}\bigl[\log\bigl\{1+\alpha(0-1 ) \bigr\}\bigr]
\\
&=&\log(1-\alpha).
\end{eqnarray*}
Thus, the set
$\mathcal{G}= \{
\log\{1+\alpha(f(\theta,F)/f(\theta_0,F_0)-1
) \}\dvtx f(\theta,F)\in\mathcal{F} \}$
has an integrable envelope.
To see this, form a sequence $(\theta_n,F_n)$ such that
\begin{eqnarray*}
g_n&\equiv&\log\bigl\{1+\alpha\bigl(f(\theta_n,F_n)/f(
\theta_0,F_0)-1 \bigr) \bigr\}
\\
&\nearrow& \sup_{\theta\in\Theta,F\in\tilde{H}} \log\bigl\{
1+\alpha\bigl(f(
\theta,F)/f(\theta_0,F_0)-1 \bigr) \bigr\} \equiv G.
\end{eqnarray*}
Then $\{g_n-\log(1-\alpha)\}_{n\in\mathbb{N}}$ is a monotone increasing
sequence of nonnegative functions.
By the monotone convergence theorem,
$P_0 g_n-\log(1-\alpha)\rightarrow P_0G -\log(1-\alpha)\leq- \log
(1-\alpha)$.
Thus we choose $G \vee-\log(1-\alpha)$ as an integrable envelope.
Also, the set $\mathcal{G}$ is Glivenko--Cantelli by a
Glivenko--Cantelli preservation theorem
\cite{MR1857319}. 

Now, by the concavity of the map $u\mapsto\log u$, and the definition
of the
WLE, we have
\begin{eqnarray*}
&&\mathbb{P}_N^\pi\log\bigl\{1+\alpha\bigl(f(\hat{
\theta}_N,\hat{F}_N)/f(\theta_0,F_0)-1
\bigr)\bigr\}
\\
&&\qquad\geq\mathbb{P}_N^\pi\bigl\{ (1-\alpha)\log(1) +\alpha
\log\bigl\{f(\hat{\theta}_N,\hat{F}_N)/f(
\theta_0,F_0)\bigr\}\bigr\}
\\
&&\qquad=\alpha\bigl\{\mathbb{P}_N^\pi\log f(\hat{
\theta}_N,\hat{F}_N)- \mathbb{P}_N^\pi
\log f(\theta_0,F_0)\bigr\} \geq0.
\end{eqnarray*}
Since $\Theta$ and $\tilde{H}$ are compact, there is a subsequence of
$(\hat{\theta}_N,\hat{F}_N)$ converging to $(\theta_\infty,F_\infty)\in
\Theta\times\tilde{H}$.
Along this subsequence, Theorem~\ref{thmfpsPiGC} implies that
\begin{eqnarray*}
0&\leq&
\mathbb{P}_N^\pi\log\bigl
\{1+\alpha\bigl(f(\hat{\theta}_N,\hat{F}_N)/f(
\theta_0,F_0)-1\bigr)\bigr\}
\\
&\rightarrow_{P^*}& 
P_{\theta_0,F_0}\bigl[\log
\bigl\{1+\alpha\bigl(f(\theta_\infty,F_\infty)/f(
\theta_0,F_0)-1 \bigr) \bigr\}\bigr]\leq0,
\end{eqnarray*}
so that $P_{\theta_0,F_0}\log\{1+\alpha(f(\theta_\infty,F_\infty)/f(\theta_0,F_0)-1 ) \}
=0$.
This is possible when $(\theta_\infty,F_\infty)=(\theta_0,F_0)$
because\vspace*{1pt}
$ (\theta, F) \mapsto P[\log\{1+\alpha(f(\theta,F)/f(\theta
_0,F_0)-1)\}]$
attains its maximum only at $(\theta_0,F_0)$.
Hence we conclude that $(\hat{\theta}_N,\hat{F}_N)$ converges to
$(\theta_0,F_0)$ in the sense of Kullback--Leibler divergence.
Since the Kullback--Leibler divergence bounds the Hellinger distance,
it follows by
Lemma A5 of~\cite{MR1463562} 
that
$d((\hat{\theta}_N,\hat{\Lambda}_N),(\theta_0,\Lambda_0))=o_{P^*}(1)$.
\end{pf}

\subsubsection{Rate of convergence}\label{sec4.2.2}
We prove the rate of convergence of the WLE is $N^{1/3}$ by
applying the rate theorem (Theorem~\ref{thmrate}) in Section~\ref{sec5}.
Since we proved the consistency of $(\hat{\theta}_N,\hat{\Lambda
}_N)$ to
$(\theta_0,\Lambda_0)$ on $S[Y]$, under Condition~\ref{condcoxi3}
we can restrict a parameter space of $\Lambda$ to
$H_M\equiv
\{\Lambda\in H\dvtx M^{-1}\leq\Lambda\leq M$, on $S[Y]\}$,
where $M$ is a positive constant such that $M^{-1}\leq\Lambda_0\leq M$
on $S[Y]$.
Define
$\mathcal{M}
\equiv\{\ell(\theta,\Lambda)\dvtx\theta\in\Theta,\Lambda\in
H_M
\}$.

%
\begin{theorem}[(Rate of convergence)]
\label{thmratecurrent}
Under Conditions~\ref{condcoxi1}--\ref{condcoxi4},
\[
d\bigl((\hat{\theta}_N,\hat{\Lambda}_N),(
\theta_0,\Lambda_0)\bigr) =O_{P^*}
\bigl(N^{-1/3}\bigr).
\]
This holds if we replace the WLE by the WLEs with adjusted
weights assuming Conditions~\ref{condest} and~\ref{condcal}.
\end{theorem}
\begin{pf}
Since the rate of convergence for the WLE is easier to verify than the
other four estimators,
we only prove the theorem for the WLE with modified calibration.
The cases for the WLEs with adjusted weights.

We proceed by verifying the conditions in Theorem~\ref{thmrate}.
Bound (\ref{eqnratecond2}) follows by Lemma~\ref{lemmaratecal} in
Section~\ref{sec5} and
Lemma A5 of~\cite{MR1463562}. 
For bound (\ref{eqnmxineqce}), we follow the proof of
(\ref{eqnmxineq}) in~\cite{MR1394975}. 
Since $\hat{\alpha}_N$ is consistent, we can specify the small neighborhood
$\mathcal{A}_{\mathrm{mc},0}$ of a zero vector such that $G_{\mathrm{mc}}(z;\alpha)$ is
contained in a small\vadjust{\goodbreak}
interval that contains 1 and consists of strictly positive numbers.
Thus, multiplying the log likelihood by a uniformly bounded quantity
$G_{\mathrm{mc}}(z;\alpha)$ only requires a slight modification of Huang's proof
of his Lemma 3.1 to obtain
$\sup_Q \log N_{[\cdot]} (\varepsilon,G\mathcal{M},L_2(Q) )
\lesssim\varepsilon^{-1}$
for $\varepsilon$ small enough where the supremum is taken over the
all discrete
probability measures and
$G\mathcal{M}
=\{G_{\mathrm{mc}}(\cdot;\alpha)\ell(\theta,\Lambda)\dvtx\alpha\in
\mathcal
{A}_{\mathrm{mc},0},\ell(\theta,\Lambda)\in\mathcal{M}\}$.
Let $G\mathcal{M}_\delta=\{m(\theta,\Lambda,\alpha)-m(\theta
_0,\Lambda
_0,\alpha)\dvtx m(\theta,\Lambda,\alpha)\in G\mathcal{M},
d ( (\theta,\Lambda), (\theta_0,\Lambda_0
) )\leq\delta\}$.
It follows\vspace*{1pt} by Lemma~3.2.2 of
\cite{MR1385671} 
that
$E^* \|\mathbb{G}_N \|_{G\mathcal{M}_\delta}
\lesssim\delta^{1/2}\{1+(\delta^{1/2}/\delta^2\sqrt{N})M\}
\equiv\phi_N(\delta)$.
Apply Theorem~\ref{thmrate} to conclude $r_N=N^{1/3}$.
\end{pf}

\subsubsection{Asymptotic normality of the estimators}\label{sec4.2.3}
We apply Theorem~\ref{thmzthm2} to derive the asymptotic distributions
of the WLEs.

%
\begin{theorem}[(Asymptotic normality)]
\label{thmancurrent} Under
Conditions~\ref{condest},~\ref{condcal},\break
\ref{condcoxi1}--\ref{condcoxi4},
\begin{eqnarray*}
\sqrt{N}(\hat{\theta}_N-\theta_0) &=&\sqrt{N}
\mathbb{P}_N^\pi\tilde{\ell}_{\theta_0,\Lambda_0}+o_{P^*}(1)
\rightsquigarrow N (0,\Sigma),
\\
\sqrt{N}(\hat{\theta}_{N,\#}-\theta_0) &=&\sqrt{N}
\mathbb{P}_N^{\pi,\#}\tilde{\ell}_{\theta_0,\Lambda
_0}+o_{P^*}(1)
\rightsquigarrow N (0,\Sigma_{\#} ), 
\end{eqnarray*}
where\vspace*{1pt}
$\#\in\{e,c,\mathrm{mc},\mathrm{cc}\}$,
$\tilde{\ell}_{\theta_0,\Lambda_0}=I^{-1}_{\theta_0,\Lambda
_0}\ell^*_{\theta_0,\Lambda_0}$
is the efficient influence function for complete data and
$\Sigma$ and $\Sigma_{\#}$ are given in Theorem~\ref{thmzthm2}.
\end{theorem}

\begin{pf}
We proceed by verifying the conditions of Theorem~\ref{thmzthm2} for
the WLE with modified calibration.
The other four cases are similar.

Condition~\ref{condwlenonreg1} is satisfied with $\beta=1/3$
by Theorems~\ref{thmconscurrent} and~\ref{thmratecurrent}.
Conditions~\ref{condwlenonreg2}--\ref{condwlenonreg4} are verified by
\cite{MR1394975} 
with
\[
\underline{h}^*(y) \equiv
\Lambda_0(y)E
\bigl[Xe^{2\theta_0^TX}O(Y|X)|Y=y\bigr]/ E\bigl[e^{2\theta
_0^TX}O(Y|X)|Y=y\bigr].
\]
Since $\mathbb{P}_N^{\pi,\mathrm{mc}}\dot{\ell}_{\hat{\theta}_{N,\mathrm{mc}},\hat
{\Lambda
}_{N,\mathrm{mc}}} =0$
by Lemma
A.5,
it remains to show that
\[
\mathbb{P}_N^{\pi,\mathrm{mc}}B_{\hat{\theta}_{N,\mathrm{mc}},\hat{\Lambda
}_{N,\mathrm{mc}}}\bigl[\underline{h}^*\bigr]=o_{P^*}\bigl(N^{-1/2}\bigr).
\]
Let $g_0\equiv\underline{h}^*\circ\Lambda_0^{-1}$ be the
composition of
$\underline{h}^*$ and the inverse of $\Lambda_0$.
Note that $\Lambda_0$ is a strictly increasing continuous function by
our assumption.
Since $g_0(\hat{\Lambda}_{N,\mathrm{mc}}(y))$ is a right continuous function
and has
exactly the same jump points as $\hat{\Lambda}_{N,\mathrm{mc}}(y)$,
by Lemma
A.5,
$\mathbb{P}_N^{\pi,\mathrm{mc}} g_0(\hat{\Lambda}_{N,\mathrm{mc}}(Y))e^{\hat{\theta
}_{N,\mathrm{mc}}^TX}
Q(Y,\Delta,X;\break\hat{\theta}_{N,\mathrm{mc}},\hat{\Lambda}_{N,\mathrm{mc}}) =0$.
By Conditions~\ref{condcoxi2}--\ref{condcoxi4}, $h^*$ has bounded derivative.
This and the assumption that $\Lambda_0$ has strictly positive
derivative by Condition~\ref{condcoxi4} imply
that $g_0$ has bounded derivative, too.
So, noting that $\underline{h}^*=g_0\circ\Lambda_0$, we have
\begin{eqnarray*}
&&
\mathbb{P}_N^{\pi,\mathrm{mc}}B_{\hat{\theta}_{N,\mathrm{mc}},\hat{\Lambda
}_{N,\mathrm{mc}}}\bigl[\underline{h}^*
\bigr] \\
&&\qquad=\mathbb{P}_N^{\pi,\mathrm{mc}} \underline{h}^*(Y)e^{\hat{\theta
}_{N,\mathrm{mc}}^TX}Q(Y,
\Delta,X;\hat{\theta}_{N,\mathrm{mc}},\hat{\Lambda}_{N,\mathrm{mc}})
\\
&&\qquad=\mathbb{P}_N^{\pi,\mathrm{mc}} \bigl\{g_0\circ
\Lambda_0(Y)-g_0\bigl(\hat{\Lambda}_{N,\mathrm{mc}}(Y)
\bigr)\bigr\} e^{\hat
{\theta
}_{N,\mathrm{mc}}^TX} Q(Y,\Delta,X;\hat{\theta}_{N,\mathrm{mc}},\hat{
\Lambda}_{N,\mathrm{mc}})
\\
&&\qquad=\bigl(\mathbb{P}_N^{\pi,\mathrm{mc}}-P_{\theta_0,\Lambda_0}\bigr) \bigl
\{g_0\circ\Lambda_0(Y)-g_0\bigl(\hat{
\Lambda}_{N,\mathrm{mc}}(Y)\bigr)\bigr\}\\
&&\qquad\quad{}\times e^{\hat
{\theta
}_{N,\mathrm{mc}}^TX} Q(Y,\Delta,X;\hat{
\theta}_{N,\mathrm{mc}},\hat{\Lambda}_{N,\mathrm{mc}})
\\
&&\qquad\quad{} + P_{\theta_0,\Lambda_0} \bigl\{g_0\circ\Lambda_0(Y)-g_0
\bigl(\hat{\Lambda}_{N,\mathrm{mc}}(Y)\bigr)\bigr\} \\
&&\hspace*{11pt}\qquad\quad{}\times e^{\hat
{\theta
}_{N,\mathrm{mc}}^TX} Q(Y,\Delta,X;
\hat{\theta}_{N,\mathrm{mc}},\hat{\Lambda}_{N,\mathrm{mc}}).
\end{eqnarray*}
Huang~\cite{MR1394975} 
showed that the second term in the display is $o_{P^*}(N^{-1/2})$.
We show that the first term in the display is also $o_{P^*}(N^{-1/2})$.
Let $C>0$ be an arbitrary constant.
Define for a fixed constant $\eta>0$,
$\mathcal{D}(\eta)
\equiv\{\psi(y,x;\theta,\Lambda)\dvtx\break d ( (\theta,\Lambda
), (\theta_0,\Lambda_0 ) )
\leq\eta,\Lambda\in H_M\}$,
where $\psi(y,\delta,x;\theta,\Lambda)\equiv
\{g_0\circ\Lambda_0(y)-\break g_0(\Lambda(y))\}\* e^{\theta^Tx}Q(y,\delta,x;\theta,\Lambda)$.
Because Huang~\cite{MR1394975} showed that $\mathcal{D}(\eta)$ is
Donsker for every
$\eta>0$ and that $\|\mathbb{G}_N\|_{\mathcal
{D}(CN^{-1/3})}=o_{P^*}(1)$,
it follows by Lemma~\ref{lemmawza5} with $\mathcal{F}_N$ replaced by
$\mathcal{D}(CN^{-1/3})$ that
$\|\mathbb{G}_N^{\pi,\mathrm{mc}}\|_{\mathcal{D}(CN^{-1/3})}=o_{P^*}(1)$.
This completes the proof.
\end{pf}
Unlike the previous example, $\ell^*_{\theta,\Lambda}$ depends on
additional unknown functions,
and the method of variance estimation used in the previous example does
not apply to the present case.
See the discussion in Section~\ref{sec6}.

\section{General results for IPW empirical processes}\label{sec5}
The IPW empirical measure and IPW empirical process inherit important
properties from the empirical measure and empirical process, respectively.
We emphasize the similarity between empirical processes and IPW
empirical processes.

\subsection{Glivenko--Cantelli theorem}\label{sec5.1}
The next theorem states that the Gli\-venko--Cantelli property for
complete data is preserved under two-phase sampling.
%
%
\begin{theorem}
\label{thmfpsPiGC}
Suppose that $\mathcal{F}$ is $P_0$-Glivenko--Cantelli.
Then
%
%
\begin{equation}
\label{eqnfpsPiGC} \bigl\|\mathbb{P}_N^\pi-P_0
\bigr\|_{\mathcal{F}}\rightarrow_{P^*} 0,
\end{equation}
where $\|\cdot\|_{\mathcal{F}}$ is the supremum norm. This also
holds if we replace $\mathbb{P}_N^\pi$ by $\mathbb{P}_N^{\pi,\#}$ with
$\#\in\{e,c,\mathrm{mc},\mathrm{cc}\}$ assuming Conditions~\ref{condest} and
\ref{condcal}.
\end{theorem}

\subsection{Rate of convergence}\label{sec5.2}
The rate of convergence of an $M$-estimator for complete data is often
established via maximal inequalities for the empirical processes.
If we follow the same line of reasoning, it is natural to derive
maximal inequalities for
IPW empirical processes, though this may require some efforts.
Fortunately, these maximal inequalities for empirical processes (or
slight modifications of them)
suffice to establish the same rate of convergence under two-phase sampling.

%
\begin{theorem}
\label{thmrate}
Let $\mathcal{M}=\{m_\theta\dvtx\theta\in\Theta\}$ be the set of
criterion functions
and define $\mathcal{M}_\delta=\{m_\theta-m_{\theta_0}\dvtx
d(\theta,\theta
_0)<\delta\}$
for some fixed $\delta>0$ where $d$ is a semimetric on the parameter
space $\Theta$.

(1)
Suppose that for every $\theta$ in a neighborhood of $\theta_0$,
%
%
\begin{equation}
\label{eqnratecond1} P_0(m_\theta-m_{\theta_0})
\lesssim-d^2(\theta,\theta_0);
\end{equation}
here $ a \lesssim b$ means $a \le K b$ for some constant $K \in
(0,\infty)$.
Assume that there exists a function $\phi_N$ such that $\delta\mapsto
\phi_N(\delta)/\delta^\alpha$
is decreasing for some $\alpha<2$ (not depending on $N$), and for
every $N$,
%
%
\begin{equation}
\label{eqnmxineq} E^*\|\mathbb{G}_N\|_{\mathcal{M}_\delta}
\lesssim\phi_N(\delta),
\end{equation}
where $\mathbb{G}_N$ is the empirical process.
If an estimator $\hat{\theta}_N$ satisfying
$\mathbb{P}_N^\pi m_{\hat{\theta}_N}\geq\mathbb{P}_N^\pi m_{\theta
_0}-O_{P^*}(r^{-2}_N)$
converges in outer probability to $\theta_0$, then $r_Nd(\hat{\theta
}_N,\theta_0)=O_{P^*}(1)$
for every sequence $r_N$ such that $r_N^2\phi_N(1/r_N)\leq\sqrt{N}$
for every $N$.

(2) Let $\#\in\{e,c,\mathrm{mc},\mathrm{cc}\}$ be fixed.
Suppose Condition~\ref{condcal} holds.
Suppose also that for every $\theta\in\Theta$ in a neighborhood of
$\theta_0$,
%
%
\begin{equation}
\label{eqnratecond2} P_0\bigl\{\tilde{G}_{\#}(V;\alpha)
(m_\theta-m_{\theta_0})\bigr\}\lesssim-d^2(\theta,
\theta_0) + |\alpha-\alpha_0|^2,
\end{equation}
where $\tilde{G}_e=\pi_0(V)/G_e$ or $\tilde{G}_\#=G_\#$ with $\#\in
\{
c,\mathrm{mc},\mathrm{cc}\}$.
Assume that
%
%
\begin{equation}
\label{eqnmxineqce} E^*\|\mathbb{G}_N\|_{\tilde{G}_\#
\mathcal
{M}_\delta}
\lesssim\phi_N(\delta),
\end{equation}
where
$\tilde{G}_\#\mathcal{M}_\delta
\equiv\{\tilde{G}_{\#}(\cdot; \alpha)f\dvtx|\alpha|\leq\delta,\alpha
\in
\mathcal{A}_N,f\in\mathcal{M}_\delta\}$
for some $\mathcal{A}_N\subset\mathcal{A}_{\#}$.
Then an estimator $\hat{\theta}_{N,\#}$ satisfying
$\mathbb{P}_N^{\pi,\#} m_{\hat{\theta}_{N,\#}}\geq\mathbb
{P}_N^{\pi,\#
} m_{\theta_0}-O_{P^*}(r^{-2}_N)$
has the same rate of convergence as $\hat{\theta}_N$ in part (1) if it
is consistent.
\end{theorem}

%
\begin{remark}
The key to establishing a general theorem for the
rate of convergence is to make use of the
boundedness of the weights in the IPW empirical process
and also deal with the dependence of the weights.
In treating independent bootstrap weights in the weighted bootstrap
(\cite{MR2202406},
Lemmas~\mbox{1--3}),
require the boundedness of bootstrap weights
because the product of an unbounded weight and a bounded function is no
longer bounded.
Our theorem exploits the boundedness of sampling indicators in the
IPW empirical processes by applying a multiplier inequality for the
case of
bounded weights (Lemma~\ref{lemmamultineq1}) to cover more general cases.
\end{remark}

The following is a multiplier inequality for bounded exchangeable weights.
Note that the sum of stochastic processes in the second term is divided
by $n^{1/2}$ rather than $k^{1/2}$.
%
%
\begin{lemma}
\label{lemmamultineq1}
For i.i.d. stochastic processes $Z_1,\ldots,Z_n$, every bounded,
exchangeable random vector\vadjust{\goodbreak} $(\xi_1,\ldots,\xi_n)$ with each $\xi_i
\in
[l,u]$ that is independent of
$Z_1,\ldots,Z_n$, and any $1\leq n_0\leq n$,
\begin{eqnarray*}
&&E \Biggl\|\frac{1}{\sqrt{n}}\sum_{i=1}^n
\xi_iZ_i \Biggr\|^*_{\mathcal{F}}
\\
&&\qquad\leq\frac{2(n_0-1)}{n}\sum_{i=1}^nE^*
\| Z_i\|_{\mathcal
{F}}E\max_{1\leq i\leq n}
\frac{\xi_i}{\sqrt{n}} +2(u-l)\max_{n_0\leq k\leq n}E \Biggl\|\frac
{1}{\sqrt{n}}\sum
_{i=n_0}^kZ_{i}
\Biggr\|^*_{\mathcal{F}}.
\end{eqnarray*}
\end{lemma}

Bound (\ref{eqnmxineqce}) is not difficult to verify in the presence of
bound
(\ref{eqnmxineq}) since $G_{\#}(\cdot;\alpha)$ is a bounded monotone
function indexed by a finite-dimensional parameter.
Bound (\ref{eqnratecond2}) may be verified through the lemma below
for some applications including the Cox model with interval censoring.
%
%
\begin{lemma}
\label{lemmaratecal}
Suppose Conditions~\ref{condest} and~\ref{condcal} hold.
Let $m_\theta$ be the log likelihood $\log p_{\theta}$ where
$p_\theta$
is the
density with dominating measure $\mu$, and $d$ is the Hellinger distance.
Then the bound (\ref{eqnratecond2}) holds.
\end{lemma}

\subsection{Donsker theorem}\label{sec5.3}
The next theorem yields weak convergence of the IPW empirical processes
under sampling without replacement.
%
%
\begin{theorem}
\label{thmfpsD}
Suppose that $\mathcal{F}$ with $\| P_0 \|_{\cal F}<\infty$
is $P_0$-Donsker and Conditions~\ref{condest} and~\ref{condcal} hold.
Then
%
%
\begin{eqnarray}
\label{eqnfpsD} \mathbb{G}_N^\pi&\rightsquigarrow&
\mathbb{G}^{\pi} \equiv\mathbb{G}+\sum_{j=1}^J
\sqrt{\nu_j}\sqrt{\frac
{1-p_j}{p_j}}\mathbb{G}_j,
\\
\mathbb{G}_N^{\pi,\#} &\rightsquigarrow&\mathbb{G}^{\pi,\#}
\equiv\mathbb{G}+\sum_{j=1}^J\sqrt{
\nu_j}\sqrt{\frac
{1-p_j}{p_j}}\mathbb{G}_j(
\cdot-Q_{\#}\cdot) 
\end{eqnarray}
in $\ell^\infty(\mathcal{F})$ where $\#\in\{e,c,\mathrm{mc},\mathrm{cc}\}$, the
$P_0$-Brownian bridge process, $\mathbb{G}$, indexed by
$\mathcal{F}$ and the $P_{0|j}$-Brownian bridge processes, $\mathbb{G}_j$,
indexed by $\mathcal{F}$ are all independent.
\end{theorem}

%
\begin{remark}
The integrability hypothesis $\| P_0\|_{\mathcal{F}}<\infty$ is
only required for the IPW empirical processes with adjusted weights.
\end{remark}

For a Donsker set $\mathcal{F}$, it follows by Theorem~\ref{thmfpsD}
and Lemma 2.3.11 of~\cite{MR1385671} 
that asymptotic equicontinuity in probability and in mean follows for
the metric that depends on the limit process.
In applications, it is of interest to have these results for the
original metric $\rho_{P_0}(f,g)=\sigma_{P_0}(f-g)$.
%
%
\begin{theorem}
\label{thmaecD}
Let $\mathcal{F}$ be Donsker and define $\mathcal{F}_\delta=\{
f-g\dvtx f,g\in\mathcal{F},\break\rho_{P_0}(f,g)<\delta\}$
for some fixed $\delta>0$.
Then, for every sequence $\delta_N\downarrow0$,
\[
E^*\bigl\|\mathbb{G}_N^\pi\bigr\|_{\mathcal{F}_{\delta
_N}}
\rightarrow0\vadjust{\goodbreak}
\]
and consequently, $\|\mathbb{G}_N^\pi\|_{\mathcal
{F}_{\delta
_N}}=o_{P^*}(1)$.
Moreover, $\|\mathbb{G}_N^{\pi,\#}\|_{\mathcal{F}_{\delta
_N}}=o_{P^*}(1)$
for $\# \in\{ e, c, \mathrm{mc}, \mathrm{cc}\}$
assuming Conditions~\ref{condest} and~\ref{condcal}.
\end{theorem}

We end this section with two important lemmas.
The first lemma is an extension of Lemma 3.3.5 of
\cite{MR1385671} 
and will be used in our proof of
Theorem~\ref{thmzthm1} to verify asymptotic equicontinuity.

%
\begin{lemma}
\label{lemmaaeccond}
Suppose ${\cal F} = \{\psi_{\theta,h}-\psi_{\theta_0,h}\dvtx\|
\theta
-\theta_0\|<\delta,h\in\mathcal{H}\}$
is $P_0$-Donsker for some $\delta>0$ and that
$\sup_{h\in\mathcal{H}}P_0(\psi_{\theta,h}-\psi_{\theta
_0,h})^2\rightarrow0$, as $\theta\rightarrow\theta_0$.
If $\hat{\theta}_N$ converges in outer probability to $\theta_0$, then
\[
\bigl\|\mathbb{G}_N^\pi(\psi_{\hat{\theta}_N,h}-
\psi_{\theta
_0,h})\bigr\|_{\mathcal{H}} = o_{P^*}(1).
\]
This also holds if we replace $\mathbb{G}_N^\pi$ by $\mathbb
{G}_N^{\pi,\#}$ with $\#\in\{e,c,\mathrm{mc},\mathrm{cc}\}$ assuming Conditions
\ref{condest} and~\ref{condcal} hold and $\| P_0 \|_{\cal F} <
\infty$.
\end{lemma}
The second lemma is used to verify asymptotic equicontinuity in the
proof of
Theorem~\ref{thmzthm2}, the first part for the IPW empirical process
and the second part for the other four IPW empirical processes with
adjusted weights.
%
%
\begin{lemma}
\label{lemmawza5}
Let $\mathcal{F}_N$ be a sequence of decreasing classes of functions
such that
$\|\mathbb{G}_N\|_{\mathcal{F}_N}=o_{P^*}(1)$.
Assume that there exists an integrable envelope for $\mathcal{F}_{N_0}$
for some $N_0$.
Then $E\|\mathbb{G}_N\|_{\mathcal{F}_N}\rightarrow0$ as
$N\rightarrow\infty$.
As a consequence, $\|\mathbb{G}_N^\pi\|_{\mathcal
{F}_N}=o_{P^*}(1)$.

Suppose, moreover, that $\mathcal{F}_N$ is $P_0$-Glivenko--Cantelli with
$\| P_0 \|_{{\cal F}_{N_1}} < \infty$ for some $N_1$,
and that every $f = f_N \in\mathcal{F}_N$ converges to zero either
pointwise or in
$L_1 (P_0)$ as $N\rightarrow\infty$.
Then $\|\mathbb{G}_N^{\pi,e}\|_{\mathcal{F}_N}=o_{P^*}(1)$,
$\|\mathbb{G}_N^{\pi,c}\|_{\mathcal{F}_N}=o_{P^*}(1)$,
$\|\mathbb{G}_N^{\pi,\mathrm{mc}}\|_{\mathcal{F}_N}=o_{P^*}(1)$ and
$\|\mathbb{G}_N^{\pi,\mathrm{cc}}\|_{\mathcal{F}_N}=o_{P^*}(1)$,
assuming Conditions~\ref{condest}\break and~\ref{condcal}.
\end{lemma}

\section{Discussion}\label{sec6}
We developed asymptotic theory for weighted likelihood
estimation under two-phase sampling, introduced and studied a new
calibration method,
centered calibration, and compared several WLE estimation methods
involving adjusted weights.
The methods of proof and general results for the IPW empirical
process are applicable to other estimation procedures.
For example, the weighted Kaplan--Meier estimator can be shown
to be asymptotically Gaussian via our Donsker theorem
(Theorem~\ref{thmfpsD}) together with the functional delta method.
A~particularly interesting application is to study asymptotic properties
of estimators that are known to be efficient under
Bernoulli sampling (e.g., estimator of~\cite{MR2125853}).
Whether or not these estimators are ``efficient'' under our sampling
scheme is an open problem; see~\cite{MR1814795} for a definition of
efficiency with non-i.i.d. data.

There are several other open problems. Variance estimation under
two-phase sampling has been restricted to the case where the asymptotic
variance is a known function up to parameters as discussed in Section
\ref{sec4}, while there are several methods available for complete data
in a general case (e.g.,~\cite{MR1693616}). In~\cite{Saegusa} the first
author has proposed and studied nonparametric bootstrap variance
estimation methods which remain valid even under model
misspecification; these results will appear elsewhere. Another
direction of research is to study (local and global) model
misspecification under two-phase sampling where missingness is by
design. An interesting open problem beyond our sampling scheme is to
study other complex survey designs. Stratified sampling without
replacement is sufficiently simple for the existing bootstrap empirical
process theory to apply. Other complex designs may provide interesting
theoretical challenges, perhaps in connection with extensions of
bootstrap empirical process theory.

\section*{Acknowledgements}

We owe thanks to Kwun Chuen Gary Chan for suggesting the modified
calibration method introduced in Section
\ref{subsubsecimprovedcalibration}. We also thank Norman Breslow for
many helpful conversations concerning two-phase sampling, and two
referees for their constructive comments and suggestions.

\begin{supplement}
\stitle{Supplementary material for ``Weighted likelihood estimation
under two-phase sampling''}
\slink[doi]{10.1214/12-AOS1073SUPP} 
\sdatatype{.pdf}
\sfilename{aos1073\_supp.pdf}
\sdescription{Due to space constraints, the proofs and technical
details have been given in the supplementary document
\cite{Saegusa-Wellner2012supp}. References here beginning with ``A.''
refer to~\cite{Saegusa-Wellner2012supp}.}
\end{supplement}


\printaddresses


\begin{thebibliography}{35}

\bibitem{MR1158522}
\begin{barticle}[mr]
\bauthor{\bsnm{Binder},~\bfnm{David~A.}\binits{D.~A.}}
(\byear{1992}).
\btitle{Fitting {C}ox's proportional hazards models from survey data}.
\bjournal{Biometrika}
\bvolume{79}
\bpages{139--147}.
\bid{doi={10.1093/biomet/79.1.139}, issn={0006-3444}, mr={1158522}}
\bptok{imsref}%
\end{barticle}
\endbibitem

\bibitem{Breslow2009a}
\begin{barticle}[author]
\bauthor{\bsnm{Breslow},~\bfnm{Norman~E.}\binits{N.~E.}},
  \bauthor{\bsnm{Lumley},~\bfnm{T}\binits{T.}},
  \bauthor{\bsnm{Ballantyne},~\bfnm{CM}\binits{C.}},
  \bauthor{\bsnm{Chambless},~\bfnm{LE}\binits{L.}} \AND
  \bauthor{\bsnm{Kulich},~\bfnm{M}\binits{M.}}
(\byear{2009}).
\btitle{Improved Horvitz--Thompson estimation of model parameters from two-phase
  stratified samples: Applications in epidemiology}.
\bjournal{Stat. Biosc.}
\bvolume{1}
\bpages{32--49}.
\bptok{imsref}%
\end{barticle}
\endbibitem

\bibitem{Breslow2009b}
\begin{barticle}[author]
\bauthor{\bsnm{Breslow},~\bfnm{Norman~E.}\binits{N.~E.}},
  \bauthor{\bsnm{Lumley},~\bfnm{T}\binits{T.}},
  \bauthor{\bsnm{Ballantyne},~\bfnm{CM}\binits{C.}},
  \bauthor{\bsnm{Chambless},~\bfnm{LE}\binits{L.}} \AND
  \bauthor{\bsnm{Kulich},~\bfnm{M}\binits{M.}}
(\byear{2009}).
\btitle{Using the whole cohort in the analysis of case-cohort data}.
\bjournal{Am. J. Epidemiol.}
\bvolume{169}
\bpages{1398--1405}.
\bptok{imsref}%
\end{barticle}
\endbibitem

\bibitem{MR2325244}
\begin{barticle}[mr]
\bauthor{\bsnm{Breslow},~\bfnm{Norman~E.}\binits{N.~E.}} \AND
  \bauthor{\bsnm{Wellner},~\bfnm{Jon~A.}\binits{J.~A.}}
(\byear{2007}).
\btitle{Weighted likelihood for semiparametric models and two-phase stratified
  samples, with application to {C}ox regression}.
\bjournal{Scand. J. Stat.}
\bvolume{34}
\bpages{86--102}.
\bid{doi={10.1111/j.1467-9469.2006.00523.x}, issn={0303-6898}, mr={2325244}}
\bptok{imsref}%
\end{barticle}
\endbibitem

\bibitem{MR2391566}
\begin{barticle}[mr]
\bauthor{\bsnm{Breslow},~\bfnm{Norman~E.}\binits{N.~E.}} \AND
  \bauthor{\bsnm{Wellner},~\bfnm{Jon~A.}\binits{J.~A.}}
(\byear{2008}).
\btitle{A {$Z$}-theorem with estimated nuisance parameters and correction note
  for: ``{W}eighted likelihood for semiparametric models and two-phase
  stratified samples, with application to {C}ox regression'' [{S}cand. {J}.
  {S}tatist. {\bf 34} (2007), no. 1, 86--102; MR2325244]}.
\bjournal{Scand. J. Stat.}
\bvolume{35}
\bpages{186--192}.
\bid{doi={10.1111/j.1467-9469.2007.00574.x}, issn={0303-6898}, mr={2391566}}
\bptok{imsref}%
\end{barticle}
\endbibitem

\bibitem{chan12}
\begin{barticle}[author]
\bauthor{\bsnm{Chan},~\bfnm{Kwun Chuen~Gary}\binits{K.~C.~G.}}
(\byear{2012}).
\btitle{Uniform improvement of empirical likelihood for missing response
  problem}.
\bjournal{Electron. J. Stat.}
\bvolume{6}
\bpages{289--302}.
\bptok{imsref}%
\end{barticle}
\endbibitem

\bibitem{MR0341758}
\begin{barticle}[mr]
\bauthor{\bsnm{Cox},~\bfnm{D.~R.}\binits{D.~R.}}
(\byear{1972}).
\btitle{Regression models and life-tables (with discussion)}.
\bjournal{J. R. Stat. Soc. Ser. B Stat. Methodol.}
\bvolume{34}
\bpages{187--220}.
\bid{issn={0035-9246}, mr={0341758}}
\bptok{imsref}%
\end{barticle}
\endbibitem

\bibitem{MR1173804}
\begin{barticle}[mr]
\bauthor{\bsnm{Deville},~\bfnm{Jean-Claude}\binits{J.-C.}} \AND
  \bauthor{\bsnm{S{\"a}rndal},~\bfnm{Carl-Erik}\binits{C.-E.}}
(\byear{1992}).
\btitle{Calibration estimators in survey sampling}.
\bjournal{J.~Amer. Statist. Assoc.}
\bvolume{87}
\bpages{376--382}.
\bid{issn={0162-1459}, mr={1173804}}
\bptok{imsref}%
\end{barticle}
\endbibitem

\bibitem{MR0053460}
\begin{barticle}[mr]
\bauthor{\bsnm{Horvitz},~\bfnm{D.~G.}\binits{D.~G.}} \AND
  \bauthor{\bsnm{Thompson},~\bfnm{D.~J.}\binits{D.~J.}}
(\byear{1952}).
\btitle{A generalization of sampling without replacement from a finite
  universe}.
\bjournal{J. Amer. Statist. Assoc.}
\bvolume{47}
\bpages{663--685}.
\bid{issn={0162-1459}, mr={0053460}}
\bptok{imsref}%
\end{barticle}
\endbibitem

\bibitem{MR1394975}
\begin{barticle}[mr]
\bauthor{\bsnm{Huang},~\bfnm{Jian}\binits{J.}}
(\byear{1996}).
\btitle{Efficient estimation for the proportional hazards model with interval
  censoring}.
\bjournal{Ann. Statist.}
\bvolume{24}
\bpages{540--568}.
\bid{doi={10.1214/aos/1032894452}, issn={0090-5364}, mr={1394975}}
\bptok{imsref}%
\end{barticle}
\endbibitem

\bibitem{MR2860827}
\begin{barticle}[mr]
\bauthor{\bsnm{Li},~\bfnm{Zhiguo}\binits{Z.}} \AND
  \bauthor{\bsnm{Nan},~\bfnm{Bin}\binits{B.}}
(\byear{2011}).
\btitle{Relative risk regression for current status data in case-cohort
  studies}.
\bjournal{Canad. J. Statist.}
\bvolume{39}
\bpages{557--577}.
\bid{doi={10.1002/cjs.10111}, issn={0319-5724}, mr={2860827}}
\bptok{imsref}%
\end{barticle}
\endbibitem

\bibitem{MR1766826}
\begin{barticle}[mr]
\bauthor{\bsnm{Lin},~\bfnm{D.~Y.}\binits{D.~Y.}}
(\byear{2000}).
\btitle{On fitting {C}ox's proportional hazards models to survey data}.
\bjournal{Biometrika}
\bvolume{87}
\bpages{37--47}.
\bid{doi={10.1093/biomet/87.1.37}, issn={0006-3444}, mr={1766826}}
\bptok{imsref}%
\end{barticle}
\endbibitem

\bibitem{LumleyRbook}
\begin{bbook}[author]
\bauthor{\bsnm{Lumley},~\bfnm{T.}\binits{T.}}
(\byear{2010}).
\btitle{Complex Surveys: A Guide to Analysis Using R}.
\bpublisher{Wiley}, \blocation{Hoboken, NJ}.
\bptok{imsref}%
\end{bbook}
\endbibitem

\bibitem{Lumley2011}
\begin{barticle}[author]
\bauthor{\bsnm{Lumley},~\bfnm{T.}\binits{T.}},
  \bauthor{\bsnm{Shaw},~\bfnm{P.~A.}\binits{P.~A.}} \AND
  \bauthor{\bsnm{Dai},~\bfnm{J.~Y.}\binits{J.~Y.}}
(\byear{2011}).
\btitle{Connections between survey calibration estimators and semiparametric
  models for incomplete data}.
\bjournal{Int. Stat. Rev.}
\bvolume{79}
\bpages{200--232}.
\bptok{imsref}%
\end{barticle}
\endbibitem

\bibitem{MR2202406}
\begin{barticle}[mr]
\bauthor{\bsnm{Ma},~\bfnm{Shuangge}\binits{S.}} \AND
  \bauthor{\bsnm{Kosorok},~\bfnm{Michael~R.}\binits{M.~R.}}
(\byear{2005}).
\btitle{Robust semiparametric {M}-estimation and the weighted bootstrap}.
\bjournal{J. Multivariate Anal.}
\bvolume{96}
\bpages{190--217}.
\bid{doi={10.1016/j.jmva.2004.09.008}, issn={0047-259X}, mr={2202406}}
\bptok{imsref}%
\end{barticle}
\endbibitem

\bibitem{MR1814795}
\begin{barticle}[mr]
\bauthor{\bsnm{McNeney},~\bfnm{Brad}\binits{B.}} \AND
  \bauthor{\bsnm{Wellner},~\bfnm{Jon~A.}\binits{J.~A.}}
(\byear{2000}).
\btitle{Application of convolution theorems in semiparametric models with
  non-i.i.d. data}.
\bjournal{J. Statist. Plann. Inference}
\bvolume{91}
\bpages{441--480}.
\bnote{Prague Workshop on Perspectives in Modern Statistical Inference:
  Parametrics, Semi-parametrics, Non-parametrics (1998)}.
\bid{doi={10.1016/S0378-3758(00)00193-2}, issn={0378-3758}, mr={1814795}}
\bptok{imsref}%
\end{barticle}
\endbibitem

\bibitem{MR1463562}
\begin{barticle}[mr]
\bauthor{\bsnm{Murphy},~\bfnm{S.~A.}\binits{S.~A.}} \AND
  \bauthor{\bparticle{van~der} \bsnm{Vaart},~\bfnm{A.~W.}\binits{A.~W.}}
(\byear{1997}).
\btitle{Semiparametric likelihood ratio inference}.
\bjournal{Ann. Statist.}
\bvolume{25}
\bpages{1471--1509}.
\bid{doi={10.1214/aos/1031594729}, issn={0090-5364}, mr={1463562}}
\bptok{imsref}%
\end{barticle}
\endbibitem

\bibitem{MR1693616}
\begin{barticle}[mr]
\bauthor{\bsnm{Murphy},~\bfnm{Susan~A.}\binits{S.~A.}} \AND
  \bauthor{\bparticle{van~der} \bsnm{Vaart},~\bfnm{Aad~W.}\binits{A.~W.}}
(\byear{1999}).
\btitle{Observed information in semi-parametric models}.
\bjournal{Bernoulli}
\bvolume{5}
\bpages{381--412}.
\bid{doi={10.2307/3318710}, issn={1350-7265}, mr={1693616}}
\bptok{imsref}%
\end{barticle}
\endbibitem

\bibitem{MR2125853}
\begin{barticle}[mr]
\bauthor{\bsnm{Nan},~\bfnm{Bin}\binits{B.}}
(\byear{2004}).
\btitle{Efficient estimation for case-cohort studies}.
\bjournal{Canad. J. Statist.}
\bvolume{32}
\bpages{403--419}.
\bid{doi={10.2307/3316024}, issn={0319-5724}, mr={2125853}}
\bptok{imsref}%
\end{barticle}
\endbibitem

\bibitem{Neyman38}
\begin{barticle}[author]
\bauthor{\bsnm{Neyman},~\bfnm{J.}\binits{J.}}
(\byear{1938}).
\btitle{Contribution to the theory of sampling human populations}.
\bjournal{J.~Amer. Statist. Assoc.}
\bvolume{33}
\bpages{101--116}.
\bptok{imsref}%
\end{barticle}
\endbibitem

\bibitem{MR1245301}
\begin{barticle}[mr]
\bauthor{\bsnm{Pr{\ae}stgaard},~\bfnm{Jens}\binits{J.}} \AND
  \bauthor{\bsnm{Wellner},~\bfnm{Jon~A.}\binits{J.~A.}}
(\byear{1993}).
\btitle{Exchangeably weighted bootstraps of the general empirical process}.
\bjournal{Ann. Probab.}
\bvolume{21}
\bpages{2053--2086}.
\bid{issn={0091-1798}, mr={1245301}}
\bptok{imsref}%
\end{barticle}
\endbibitem

\bibitem{Prentice86}
\begin{barticle}[author]
\bauthor{\bsnm{Prentice},~\bfnm{R.~L.}\binits{R.~L.}}
(\byear{1986}).
\btitle{A case-cohort design for epidemiologic cohort studies and disease
  prevention trials}.
\bjournal{Biometrika}
\bvolume{73}
\bpages{1--11}.
\bptok{imsref}%
\end{barticle}
\endbibitem

\bibitem{MR1294730}
\begin{barticle}[mr]
\bauthor{\bsnm{Robins},~\bfnm{James~M.}\binits{J.~M.}},
  \bauthor{\bsnm{Rotnitzky},~\bfnm{Andrea}\binits{A.}} \AND
  \bauthor{\bsnm{Zhao},~\bfnm{Lue~Ping}\binits{L.~P.}}
(\byear{1994}).
\btitle{Estimation of regression coefficients when some regressors are not
  always observed}.
\bjournal{J. Amer. Statist. Assoc.}
\bvolume{89}
\bpages{846--866}.
\bid{issn={0162-1459}, mr={1294730}}
\bptok{imsref}%
\end{barticle}
\endbibitem

\bibitem{Saegusa}
\begin{bmisc}[author]
\bauthor{\bsnm{Saegusa},~\bfnm{Takumi}\binits{T.}}
(\byear{2012}).
\bhowpublished{Weighted likelihood estimation under two-phase sampling.
Ph.D. thesis,
Univ. Washington, Seattle, WA}.
\bptok{imsref}%
\end{bmisc}
\endbibitem

\bibitem{Saegusa-Wellner2012supp}
\begin{bmisc}[author]
\bauthor{\bsnm{Saegusa},~\bfnm{Takumi}\binits{T.}} \AND
  \bauthor{\bsnm{Wellner},~\bfnm{Jon~A.}\binits{J.~A.}}
(\byear{2012}).
\bhowpublished{Supplement to ``Weighted likelihood estimation under
  two-phase sampling.'' DOI:\doiurl{10.1214/12-AOS1073SUPP}}.
\bptok{imsref}%
\end{bmisc}
\endbibitem

\bibitem{Saegusa-Wellner2012TechRep}
\begin{bmisc}[author]
\bauthor{\bsnm{Saegusa},~\bfnm{Takumi}\binits{T.}} \AND
  \bauthor{\bsnm{Wellner},~\bfnm{Jon~A.}\binits{J.~A.}}
(\byear{2012}).
\bhowpublished{Weighted likelihood estimation under two-phase sampling.
Technical Report 592,
Dept. Statistics, Univ. Washington, Seattle, WA.
Available at arXiv:\arxivurl{1112.4951}}.
\bptok{imsref}%
\end{bmisc}
\endbibitem

\bibitem{MR924857}
\begin{barticle}[mr]
\bauthor{\bsnm{Self},~\bfnm{Steven~G.}\binits{S.~G.}} \AND
  \bauthor{\bsnm{Prentice},~\bfnm{Ross~L.}\binits{R.~L.}}
(\byear{1988}).
\btitle{Asymptotic distribution theory and efficiency results for case-cohort
  studies}.
\bjournal{Ann. Statist.}
\bvolume{16}
\bpages{64--81}.
\bid{doi={10.1214/aos/1176350691}, issn={0090-5364}, mr={0924857}}
\bptok{imsref}%
\end{barticle}
\endbibitem

\bibitem{MR2836413}
\begin{barticle}[mr]
\bauthor{\bsnm{Tan},~\bfnm{Z.}\binits{Z.}}
(\byear{2011}).
\btitle{Efficient restricted estimators for conditional mean models with
  missing data}.
\bjournal{Biometrika}
\bvolume{98}
\bpages{663--684}.
\bid{doi={10.1093/biomet/asr007}, issn={0006-3444}, mr={2836413}}
\bptok{imsref}%
\end{barticle}
\endbibitem

\bibitem{MR1915446}
\begin{bincollection}[mr]
\bauthor{\bparticle{van~der} \bsnm{Vaart},~\bfnm{Aad}\binits{A.}}
(\byear{2002}).
\btitle{Semiparametric statistics}.
In \bbooktitle{Lectures on Probability Theory and Statistics ({S}aint-{F}lour,
  1999)}.
\bseries{Lecture Notes in Math.}
\bvolume{1781}
\bpages{331--457}.
\bpublisher{Springer}, \blocation{Berlin}.
\bid{mr={1915446}}
\bptok{imsref}%
\end{bincollection}
\endbibitem

\bibitem{MR1857319}
\begin{bincollection}[mr]
\bauthor{\bparticle{van~der} \bsnm{Vaart},~\bfnm{Aad}\binits{A.}} \AND
  \bauthor{\bsnm{Wellner},~\bfnm{Jon~A.}\binits{J.~A.}}
(\byear{2000}).
\btitle{Preservation theorems for {G}livenko--{C}antelli and uniform
  {G}livenko--{C}antelli classes}.
In \bbooktitle{High Dimensional Probability, {II} ({S}eattle, {WA}, 1999)}.
\bseries{Progress in Probability}
\bvolume{47}
\bpages{115--133}.
\bpublisher{Birkh\"auser}, \blocation{Boston, MA}.
\bid{mr={1857319}}
\bptok{imsref}%
\end{bincollection}
\endbibitem

\bibitem{MR1652247}
\begin{bbook}[mr]
\bauthor{\bparticle{van~der} \bsnm{Vaart},~\bfnm{A.~W.}\binits{A.~W.}}
(\byear{1998}).
\btitle{Asymptotic Statistics}.
\bseries{Cambridge Series in Statistical and Probabilistic Mathematics}
\bvolume{3}.
\bpublisher{Cambridge Univ. Press}, \blocation{Cambridge}.
\bid{mr={1652247}}
\bptok{imsref}%
\end{bbook}
\endbibitem

\bibitem{MR1385671}
\begin{bbook}[mr]
\bauthor{\bparticle{van~der} \bsnm{Vaart},~\bfnm{Aad~W.}\binits{A.~W.}} \AND
  \bauthor{\bsnm{Wellner},~\bfnm{Jon~A.}\binits{J.~A.}}
(\byear{1996}).
\btitle{Weak Convergence and Empirical Processes: With Applications to Statistics}.
\bpublisher{Springer}, \blocation{New York}.
\bid{mr={1385671}}
\bptok{imsref}%
\end{bbook}
\endbibitem

\bibitem{White1982}
\begin{barticle}[author]
\bauthor{\bsnm{White},~\bfnm{J.~Emily}\binits{J.~E.}}
(\byear{1986}).
\btitle{A two stage design for the study of the relationship between a rare
  exposure and and a rare disease}.
\bjournal{Am. J. Epidemiol.}
\bvolume{115}
\bpages{119--128}.
\bptok{imsref}%
\end{barticle}
\endbibitem

\bibitem{Zheng-Little04}
\begin{barticle}[author]
\bauthor{\bsnm{Zheng},~\bfnm{Hui}\binits{H.}} \AND
  \bauthor{\bsnm{Little},~\bfnm{Roderick J.~A.}\binits{R.~J.~A.}}
(\byear{2004}).
\btitle{Penalized spline nonparametric mixed models for inference about a
  finite population mean from two-stage samples}.
\bjournal{Survey Methodology}
\bvolume{30}
\bpages{209--218}.
\bptok{imsref}%
\end{barticle}
\endbibitem

\end{thebibliography}
\end{document}